\RequirePackage{fix-cm}
\documentclass[a4paper, 12pt, leqno]{article}

\usepackage{latexsym, amsmath, amsthm, amssymb, graphicx, graphics, epsfig}
\usepackage{bm}
\usepackage{tikz-cd}
\usepackage{fullpage}

\newcommand{\comment}[1]{}

\usepackage[T1]{fontenc}
\usepackage{latexsym}
\usepackage{manfnt}
\usepackage{palatino}
\usepackage[sc, osf]{mathpazo} 
\usepackage{eulervm}
\usepackage{euscript}

\newcommand{\sub}[1]{\vspace{8pt}\textbf{#1}\hspace*{0.5em}}
\newcommand{\upskip}{\vspace{-7pt}}
\newcommand{\uppskip}{\vspace{-14pt}}

\newtheorem{Thm}{Theorem}

\newtheorem{thm}{Theorem}[section]
\newtheorem*{thm*}{Theorem}
\newtheorem{lem}[thm]{Lemma}
\newtheorem*{lem*}{Lemma}

\newtheorem*{cor*}{Corollary}
\newtheorem{prop}[thm]{Proposition}
\newtheorem*{prop*}{Proposition}
\newtheorem*{claim*}{Claim}

\theoremstyle{definition}

\newtheorem*{rem}{Remark}

\newtheorem*{Conv*}{Convention}

\numberwithin{equation}{section}

\renewcommand{\proof}{\vspace{-5pt}\noindent\textit{\textbf{Proof. }}}
\newcommand{\prooff}[1]{\vspace{-5pt}\noindent\textit{\textbf{Proof{#1}. }}}
\renewcommand{\endproof}{$\square$}

\newcommand{\q}[1]{{``#1''}}
\newcommand{\br}[1]{\left(#1\right)}

\newcommand{\bra}[1]{\left\langle #1 \right\rangle}
\newcommand{\set}[1]{\left\{#1\right\}}
\newcommand{\norm}[1]{\left\|#1\right\|}
\newcommand{\case}[1]{\begin{cases}#1\end{cases}}
\newcommand{\inv}{^{-1}}

\newcommand{\restr}[2]{\left. #1 \right|_{#2}}

\newcommand{\Em}{E^-}
\newcommand{\Ep}{E^+}

\newcommand{\MS}{M\times S^1}
\newcommand{\pM}{\p M}
\newcommand{\pMS}{\pM\times S^1}
\newcommand{\EY}{E_{\p}}
\newcommand{\EmY}{\Em_{\p}}
\newcommand{\EpY}{\Ep_{\p}}

\newcommand{\Eh}{\widehat{E}}
\newcommand{\Fh}{\widehat{F}}

\newcommand{\EEm}{\E^-}

\newcommand{\EEY}{\E_{\p}}
\newcommand{\EEmY}{\EEm_{\p}}

\newcommand{\cut}{^{\mathrm{cut}}}

\newcommand{\matr}[1]{
\begin{pmatrix}
#1
\end{pmatrix}}

\newcommand{\smatr}[1]{\br{ \begin{smallmatrix}#1\end{smallmatrix} } }

\newcommand{\R}{\mathbb R}
\newcommand{\N}{\mathbb N}
\newcommand{\Z}{\mathbb Z}
\newcommand{\CC}{\mathbb C}

\newcommand{\BB}{\mathcal B}
\newcommand{\Cl}{\mathcal C}

\newcommand{\Reg}{\mathcal R}

\newcommand{\FR}{\mathcal{FR}}

\newcommand{\B}{\mathcal B}

\newcommand{\F}{\mathcal F}

\newcommand{\E}{\mathcal E}
\newcommand{\A}{\mathcal A}

\renewcommand{\L}{\mathcal L}

\newcommand{\D}{\mathcal D}

\newcommand{\overbar}[1]{\overline{#1\mkern-1mu}\mkern 1mu}
\DeclareMathOperator{\Ell}{Ell}
\newcommand{\Ellt}{\overbar{\Ell}}
\newcommand{\Ellp}{\Ellt^{+}}
\newcommand{\Ellm}{\Ellt^{-}}
\newcommand{\Ellz}{\widetilde{\Ell}}

\DeclareMathOperator{\Dir}{Dir}
\DeclareMathOperator{\Dirt}{\overbar{Dir}}

\newcommand{\Dirp}{\overbar{\Dir}^{+}}
\newcommand{\Dirm}{\overbar{\Dir}^{-}}

\DeclareMathOperator{\Opd}{Op^d}
\DeclareMathOperator{\Op1}{Op^1}
\newcommand{\Opt}{\widetilde{Op}}
\DeclareMathOperator{\Proj}{Proj}
\DeclareMathOperator{\im}{Im}
\DeclareMathOperator{\Ker}{Ker}

\DeclareMathOperator{\Gr}{Gr}
\DeclareMathOperator{\Grr}{Gr^{(2)}}

\DeclareMathOperator{\Hom}{Hom}
\DeclareMathOperator{\End}{End}

\DeclareMathOperator{\Id}{Id}

\DeclareMathOperator{\supp}{supp}

\DeclareMathOperator{\rank}{rank}
\DeclareMathOperator{\dom}{dom}

\DeclareMathOperator{\dist}{dist}

\DeclareMathOperator{\spf}{\mathsf{sf}}
\renewcommand{\sp}[1]{\spf\,(#1)}

\renewcommand{\d}{d}

\newcommand{\p}{\partial}

\newcommand{\Cinf}{C^{\infty}}

\newcommand{\sa}{^{\mathrm{sa}}}
\newcommand{\ort}{^{\mathrm{ort}}}
\newcommand{\Pmort}{P^-_{\mathrm{ort}}}
\newcommand{\ann}{_{\mathrm{ann}}}

\renewcommand{\leq}{\leqslant}
\renewcommand{\geq}{\geqslant}

\newcommand{\hookto}{\hookrightarrow}

\title{Self-adjoint local boundary problems on compact surfaces. I. Spectral flow} 
\author{Marina\,Prokhorova}
\date{}

\begin{document}

\maketitle

\upskip
\begin{abstract}
The paper deals with first order self-adjoint elliptic differential operators
on a compact oriented surface with non-empty boundary.
We consider such operators with self-adjoint local boundary conditions.
The paper is focused on paths in the space of such operators 
connecting two operators conjugated by a unitary automorphism.
We compute the spectral flow for such paths
in terms of the topological data over the boundary.
The second result is the universality of the spectral flow:
we show that the spectral flow is a universal additive invariant for such paths,
if the vanishing on paths of invertible operators is required.
\end{abstract}

\let\oldnumberline=\numberline
\renewcommand{\numberline}{\vspace{-16pt}\oldnumberline}
\addtocontents{toc}{\protect\renewcommand{\bfseries}{}}

\upskip
\tableofcontents
\bigskip

\clearpage

\section{Introduction}

\upskip
\sub{Conventions.}
Throughout the paper a \q{Hilbert space} always means a separable complex Hilbert space of infinite dimension
and a \q{surface} always means a smooth compact oriented surface with non-empty boundary.
By the \q{symbol of a differential operator} we always mean its principal symbol.

\upskip
\sub{The spectral flow for unbounded operators.}
A closed linear operator $A$ on a Hilbert space $H$ is a (not necessarily bounded) linear operator 
acting from a linear subspace $\dom(A)\subset H$ to $H$ such that its graph is closed in $H\oplus H$.
The natural topology on the space of closed operators on $H$ is the so-called graph topology
induced by the metric $\delta\br{A_1, A_2} = \norm{P_1-P_2}$,
where $P_i$ denotes the orthogonal projection of $H\oplus H$ onto the graph of $A_i$. 

The space $\FR\sa(H)$ of regular (that is, closed and densely defined) 
Fredholm self-adjoint operators on $H$ equipped with the graph topology
is path-connected, and its fundamental group is isomorphic to $\Z$. 
This isomorphism is given by the 1-cocycle on $\FR\sa(H)$ called the spectral flow.
Roughly speaking, the spectral flow counts with signs 
the number of eigenvalues passing through zero from the start of the path to its end
(the eigenvalues passing from negative values to positive one are counted with a plus sign,
and the eigenvalues passing in the other direction are counted with a minus sign).
See \cite{BLP-04} for a rigorous definition.

\sub{Local boundary value problems on a compact manifold.}
Let $M$ be a smooth compact Riemannian manifold with non-empty boundary $\pM$,
and let $E$ be a Hermitian vector bundle over $M$.
Denote by $\EY$ the restriction of $E$ to $\pM$.
Let $A$ be a first order formally self-adjoint elliptic differential operator acting on sections of $E$.
We consider only local, or classical, boundary conditions for $A$, 
which are defined by smooth subbundles of $\EY$
(in particular, boundary conditions defined by spectral projections are not allowed).
A smooth subbundle $L$ of $\EY$ defines the unbounded operator $A_L$ 
on the space $L^2(E)$ of square-integrable sections of $E$; its domain is
\begin{equation*}
  \dom\br{A_L} = \set{ u\in H^1(E)\colon \restr{u}{\pM} \mbox{ is a section of } L },
\end{equation*}
where $H^1(E)$ denotes the first order Sobolev space of sections of $E$.
(More precisely, the boundary condition above means that the boundary trace of $u$ is an element of $H^{1/2}(L)$,
see explanation in Section \ref{sec:bound_pr}.)

Let $n$ denote the outward conormal to the boundary $\pM$.
The conormal symbol $\sigma(n)$ of $A$ is self-adjoint 
and thus defines a symplectic structure on the fibers of $\EY$
given by the symplectic 2-form $\omega_x\colon E_x\otimes E_x\to\CC$,
$\omega_x(u,v) = \bra{i\sigma(n_x)\, u, \, v}$ for $x\in \pM$.
Green's formula for $A$ can be written as 
\[
\bra{Au,v}_{L^2(E)} - \bra{u,Av}_{L^2(E)} = \bra{i\sigma(n)\restr{u}{\pM},\restr{v}{\pM}}_{L^2(\EY)}
 = \omega\br{\restr{u}{\pM},\restr{v}{\pM}} 
\text{ for } u,v\in H^1(E).
\]

Since $A$ is elliptic, $\sigma(tn_x+\xi) = t\sigma(n_x) + \sigma(\xi)$ is non-degenerate
for every $t\in\R$, $x\in \pM$, and non-zero cotangent vector $\xi\in T^*_x\pM$
(loosely speaking, we identify here $T^*_x\pM$ with a subspace of $T^*_x M$; 
see Section \ref{sec:bound_pr} for more precise statements).
Hence the fiber endomorphism (the symbol of a tangential operator) $b(\xi) = \sigma(n_x)^{-1} \sigma(\xi)$ 
has no eigenvalues on the real axis. 
The generalized eigenspaces $\Ep(\xi)$ and $\Em(\xi)$ of $b(\xi)$ 
corresponding to eigenvalues with positive, resp. negative imaginary part
are Lagrangian subspaces of $E_x$ 
(that is $\sigma(n_x)$ takes them to their orthogonal complements).

A local boundary condition $L$ is called an elliptic, 
or Shapiro-Lopatinskii, boundary condition for $A$ 
if for every $x\in \pM$ the subspace $L_x$ is complementary for each $\Ep(\xi)$, that is
\begin{equation}\label{eq:L-ell0}
	L_x \cap \Ep(\xi) = 0 \;\text{ and }\; L_x+\Ep(\xi) = E_x \;
	\text{ for every non-zero } \xi\in T_x^*\pM. 
\end{equation}
If $L$ is elliptic for $A$, then 
$A_L$ is a closed operator on $H=L^2(E)$ with compact resolvents. 
If, in addition, $L$ is a Lagrangian subbundle of $\EY$, that is
\begin{equation}\label{eq:Lsa0}
\sigma(n) L = L^{\bot},
\end{equation}
then $A_L$ is self-adjoint.

We denote by $\Ellt(E)$ the space of all such pairs $(A,L)$ 
equipped with the $C^{1}$-topology on symbols  of operators, the $C^{0}$-topology on their free terms, 
and the $C^{1}$-topology on boundary conditions.
The natural inclusion $\Ellt(E)\hookto\FR\sa(H)$ taking $(A,L)$ to $A_L$
is continuous, see Proposition \ref{prop:AL2}.
Thus the spectral flow $\sp{\gamma}$ is defined for every continuous path $\gamma\colon[0,1]\to\Ellt(E)$. 

We are interested in computing the spectral flow for paths in $\Ellt(E)$ with conjugate ends, 
since the spectral flow is homotopy invariant for such paths.
In the paper we consider the simplest non-trivial case, 
namely the case of a \textbf{two-dimensional} manifold $M$.
For such $M$, we compute the spectral flow in terms of the topological data 
extracted from the corresponding one-parameter family of operators and boundary conditions.

\sub{Local boundary value problems on a surface.}
Let $M$ be an oriented smooth compact surface. 
Then $T^*_x\pM \setminus\set{0}$ consists of only two rays,
so $\EY$ can be naturally decomposed into the direct sum $\EpY\oplus\EmY$ of two Lagrangian subbundles.
Their fibers can be written as $\Ep_x = \Ep(\xi)$ and $\Em_x = \Em(\xi)$,
where $(n_x,\xi)$ is a positive oriented frame in $T^{*}_x M$.

The identity $\Ep(\xi) = \Em(-\xi)$ together with \eqref{eq:Lsa0}
allows to simplify the ellipticity condition \eqref{eq:L-ell0}.
Namely, a self-adjoint elliptic local boundary condition for $A$
is a Lagrangian subbundle $L$ of $\EY$ satisfying 
\begin{equation*}
	L \cap \EpY = L \cap \EmY = 0.
\end{equation*}
We show in Proposition~\ref{prop:T} that such subbundles $L$
are in a one-to-one correspondence with self-adjoint bundle automorphisms $T$ of $\EmY$.
This correspondence is given by the rule
\begin{equation}\label{eq:PT0}
	      L=\Ker P_T \; \mbox{ with } \; P_T = P^+\br{ 1+i\sigma(n)\inv TP^- },
\end{equation}
where $P^+$ denotes the bundle projection of $\EY$ onto $\EpY$ along $\EmY$ and $P^- = 1-P^+$.
Moreover, $P_T$ is the bundle projection of $\EY$ onto $\EpY$ along $L$,
so a local boundary condition given by $L$ can be written equivalently in the form $P_T(\restr{u}{\pM}) = 0$
using a bundle projection (a particular case of a pseudo-differential projection), 
as is customary in the study of boundary value problems.

If $A$ is the Dirac operator, then $\EpY$ and $\EmY$ are mutually orthogonal; 
in this case $L$ can be written (fiber-wise) as
$L=\set{ u^+\oplus u^-\in \EpY\oplus\EmY \colon \; i\sigma(n)u^+=Tu^- }$.

\sub{The topological data.}
We associate with an element $(A,L)\in\Ellt(E)$ the vector subbundle $F=F(A,L)$ of $\EmY$,
whose fibers $F_x$, $x\in\pM$, are spanned by the generalized eigenspaces of $T_x$ corresponding to negative eigenvalues.

Let $\gamma\colon[0,1]\to\Ellt(E)$, $\gamma = (A_t,L_t)$ 
be a continuous path such that $\gamma(1)=g\gamma(0)$ for some smooth unitary bundle automorphism $g$ of $E$.
With every such pair $(\gamma,g)$ we associate the vector bundle $\F(\gamma,g)$ over the product $\pMS$ as follows.
The one-parameter family $(F_t)$ of subbundles $F_t=F(A_t,L_t)$ of $\EY$ 
defines the vector bundle over $\pM\times[0,1]$.
The condition $\gamma(1)=g\gamma(0)$ implies $F_1=gF_0$.
Gluing $F_1$ with $F_0$ twisted by $g$, 
that is identifying $(u,1)$ with $(gu,0)$ for every $u\in F_0$,
we obtain the vector bundle $\F=\F(\gamma,g)$ over $\pMS$.

The product $\pMS$ is a disjoint union of tori.
The orientation on $M$ induces the orientation on $\pMS$.
Evaluating the first Chern class of the vector bundle $\F(\gamma,g)$ on the fundamental class of $\pMS$, 
we obtain the integer-valued invariant 
\begin{equation*}
 \Psi(\gamma,g) = c_1(\F(\gamma,g))[\pMS] = \sum_{j=1}^{m}c_1(\F_j)[\pM_j\times S^1],
\end{equation*}
where $\pM_j$, $j=1\ldots,m$, are the boundary components and $\F_j$ is the restriction of $\F$ to $\pM_j$.

\sub{The spectral flow formula.} 
The first main result of the paper is the following formula.

\begin{Thm}[Theorem \ref{thm:sf}]\label{Thm:sf}
Let $\gamma\colon[0,1]\to\Ellt(E)$ be a continuous path such that 
$\gamma(1)=g\gamma(0)$ for some smooth unitary bundle automorphism $g$ of $E$.
Then the spectral flow for $\gamma$ can be computed in terms of the topological data over the boundary:
\begin{equation*}
 \sp{\gamma} = c_1(\F(\gamma,g))[\pMS] = \Psi(\gamma,g),
\end{equation*}
\end{Thm}

\upskip 
This result was first announced by the author in \cite[Section 8]{Pr}
(up to multiplication by an integer constant depending only on the homotopy type of $M$)
and then in \cite{Pr14}.

Note that we \textit{do not} require the weak inner unique continuation property for the operators $\gamma(t)$.
While Dirac operators always have this property,
for general first order self-adjoint elliptic operators this is not necessarily so.

\sub{Universality of the spectral flow.} 
The second main result of the paper is universality of the spectral flow 
for paths in $\Ellt(E)$ with conjugate ends. 

Denote by $U(E)$ the group of smooth unitary bundle automorphisms of $E$.
For $g\in U(E)$ we denote by $\Omega_g\Ellt(E)$ the space of continuous paths 
$\gamma\colon[0,1]\to\Ellt(E)$ such that $\gamma(1)=g\gamma(0)$, equipped with the compact-open topology.

Recall that every complex vector bundle over $M$ is trivial and that $\Ellt(E)$ is empty for bundles $E$ of odd rank.
Denote by $2k_M$ the trivial vector bundle of rank $2k$ over $M$.

\begin{Thm}[Theorem \ref{thm:uni-sf}]\label{Thm:uni-sf}
Let $\Lambda$ be a commutative monoid.
Suppose that we associate an element $\Phi(\gamma,g)\in\Lambda$
with every path $\gamma\in\Omega_g\Ellt(2k_M)$ for every $k\in\N$ and every $g\in U(2k_M)$. 
Then the following two conditions are equivalent:
\begin{enumerate}\upskip
	\item $\Phi$ is homotopy invariant, additive with respect to direct sums, 
and vanishing on paths of invertible operators. 
	\item $\Phi$ has the form $\Phi(\gamma,g)  = \spf(\gamma)\cdot \lambda$
for some (invertible) constant $\lambda\in\Lambda$.
\end{enumerate}
\end{Thm}

\upskip 
The homotopy invariance here is understood as the invariance with respect to a change
of a path in the space $\Omega_g \Ellt(E)$ of all paths in $\Ellt(E)$ 
with ends conjugated by a \textit{fixed} unitary automorphism $g$ of $E$.
In other words, $\Phi$ is constant on path connected components of $\Omega_g\Ellt(2k_M)$.

By vanishing on paths of invertible operators we mean that 
$\Phi$ vanishes on $\Omega_g\Ellt^0(2k_M)$ for every $k$ and $g$,
where $\Ellt^0(E)$ denotes the subspace of $\Ellt(E)$ consisting of all pairs $(A,L)$
such that the unbounded operator $A_L$ is invertible 
(or, what is the same, has no zero eigenvalues). 
 
A similar result holds also for invariants $\Phi$ defined only on loops $\gamma\in\Omega\Ellt(2k_M)$,
see previous preprint \cite[Theorem 2]{Pr17}.
We exclude this result from the current version of the paper to make the exposition more clear.
The generalization of that result can be found in the next paper; see \cite[Theorem 11.5]{Pr18}.

It is known that the spectral flow is a universal homotopy invariant for loops in the space $\FR\sa(H)$,
and that the spectral flow is additive with respect to direct sums and vanishes on loops of invertible operators.
But the space $\Ellt(E)$ is only tiny part of $\FR\sa(L^2(E))$.
Universality is usually lost after passing to a subspace,
so we cannot expect the spectral flow to be a universal invariant for loops in $\Ellt(E)$.
Indeed, for any given $E$ the map $\spf\colon[S^1,\Ellt(E)]\to\Z$ is not injective.
It is surprising that universality can still be restored by considering all vector bundles over $M$ together.

\sub{Universality of $\Psi$.}
The proofs of Theorems A and B are based on the following result, 
which we prove in Section \ref{sec:uni-Psi} using topological means only.

Denote by $\Ellp(E)$, resp. $\Ellm(E)$ the subspace of $\Ellt(E)$
consisting of all $(A,L)$ with positive, resp. negative definite $T$,
where the bundle automorphism $T=T(A,L)$ is defined by formula \eqref{eq:PT0}.

\begin{Thm}[Theorem \ref{thm:uni-Psi}]\label{Thm:uni-Psi}
Let $\Lambda$ be a commutative monoid.
Suppose that we associate an element $\Phi(\gamma,g)\in\Lambda$
with every path $\gamma\in\Omega_g\Ellt(2k_M)$ 
for every $k\in\N$ and $g\in U(2k_M)$. 
Then the following two conditions are equivalent:
\begin{enumerate}\upskip
	\item $\Phi$ is homotopy invariant, additive with respect to direct sums, 
and vanishing on constant loops in $\Ellt(2k_M)$ 
and on paths in $\Ellp(2k_M)$, $\Ellm(2k_M)$ for every $k$.
	\item $\Phi$ has the form $\Phi(\gamma,g)  = \Psi(\gamma,g)\cdot \lambda$
for some (invertible) constant $\lambda\in\Lambda$.
\end{enumerate}
\end{Thm}

\upskip 
The direction ($2\Rightarrow 1$) follows immediately from the properties of $\Psi$.
To prove ($1\Rightarrow 2$), we first notice that 
if an additive homotopy invariant $\Phi$
vanishes on $\Omega_g\Ellp(2k_M)$ and $\Omega\Ellm(2k_M)$ for every $k$ and $g$,
then it depends only on the class of $\F(\gamma,g)$ in $K^0(\pM\times S^1)$.
Next we show that vanishing of $\Phi$ on $\Omega_g\Ellm(2k_M)$ cancels the image $G^{\p}$ 
of the homomorphism $K^0(\MS) \to K^0(\pMS)$ induced by the \mbox{embedding} \mbox{$\pMS\hookto\MS$}.
Similarly, vanishing of $\Phi$ on constant loops cancels the image $G^*$ 
of the homomorphism $K^0(\pM) \to K^0(\pMS)$ induced by the projection $\pMS\to\pM$.
The subgroup of $K^0(\pM\times S^1)$ spanned by $G^{\p}$ and $G^*$ 
is the kernel of the surjective homomorphism $\psi\colon K^0(\pMS)\to\Z$,
which is given by the rule $\psi[V] = c_1(V)[\pMS]$ for every vector bundle $V$ over $\pMS$.
It follows that $\Phi$ factors through $\psi$,
that is $\Phi(\gamma,g) = \vartheta\circ\psi[\F(\gamma,g)] = \vartheta(\Psi(\gamma,g))$ 
for some monoid homomorphism $\vartheta\colon\Z\to\Lambda$.
It remains to take $\lambda=\vartheta(1)$.

\sub{Invertible operators.}
Obviously, every constant loop $\gamma\in\Omega\Ellt(E)$ 
is homotopic to a constant loop $\gamma'\in\Omega\Ellt^0(E)$: 
one can just add a constant to the correspondent operator.

Denote by $\Dirt(E)$ the subspace of $\Ellt(E)$ consisting of all pairs $(A,L)$ such that 
$A$ is a Dirac operator which is odd with respect to the chiral decomposition.
Two subspaces of $\Dirt(E)$ play special role:
$\Dirp(E) = \Dirt(E)\cap\Ellp(E)$ and $\Dirm(E) = \Dirt(E)\cap\Ellm(E)$.
It can be easily seen that the unbounded operator $A_L$ is invertible for every 
$(A,L)\in\Dirp(E)$ or $\Dirm(E)$, see Proposition \ref{prop:zero1} for detail.
Thus both $\Omega_g\Dirp(E)$ and $\Omega_g\Dirm(E)$ are subspaces of $\Omega_g\Ellt^0(E)$.

\sub{Deformation retraction.}
We show in Section \ref{sec:retr} that the natural embedding 
$\Dirt(E)\hookto\Ellt(E)$ is a homotopy equivalence. 
Moreover, we construct a deformation retraction of $\Ellt(E)$ onto a subspace of $\Dirt(E)$
preserving the vector bundles $\EmY(A)$ and $F(A,L)$, see Proposition \ref{prop:retr1t}. 
Similarly, we construct a deformation retraction of $\Omega_g\Ellt(E)$ onto a subspace of $\Omega_g\Dirt(E)$
preserving the vector bundles $\EEmY(\gamma,g)$ and $\F(\gamma,g)$,
see Proposition \ref{prop:retr2}.
Restricting the last retraction to the special subspaces defined above,
we obtain a deformation retraction of $\Omega_g\Ellp(E)$ onto a subspace of $\Omega_g\Dirp(E)$
and a deformation retraction of $\Omega_g\Ellm(E)$ onto a subspace of $\Omega_g\Dirm(E)$.

In particular, every path connected component of $\Omega_g\Ellp(E)$, resp. $\Omega_g\Ellm(E)$
contains an element of $\Omega_g\Dirp(E)$, resp. $\Omega_g\Dirm(E)$.
It follows that every function $\Phi$ satisfying the first condition of Theorem B
should satisfy also the first condition of Theorem C.
We use this result to deduce Theorems A and B from Theorem C.

\sub{Proof of Theorem A.}
To prove Theorem A, we use the homotopy invariance of the spectral flow,
its additivity with respect to direct sums,
and vanishing of the spectral flow on paths of invertible operators.
In other words, the spectral flow considered as a function $\spf\colon\Omega_g\Ellt(E)\to\Z$ 
satisfies the first condition of Theorem B and thus of Theorem C, with $\Phi=\spf$ and $\Lambda=\Z$.
Theorem C implies that there is an integer constant $\lambda\in\Z$ depending only on $M$
such that $\spf(\gamma)  = \Psi(\gamma,g)\cdot \lambda$ 
for every $\gamma\in\Omega_g\Ellt(E)$. 

It remains to find the factor $\lambda = \lambda_M$.
Simple reasoning shows that $\lambda_M$ depends only on the diffeomorphism type of $M$.
We then reduce the computation of $\lambda_M$ to the case of the annulus,
see Lemma \ref{lem:lambda}.
The factor $\lambda\ann$ was computed by the author in \cite{Pr} by direct evaluation.
This gives $\lambda_M=\lambda\ann=1$ for any surface $M$
and completes the proof of Theorem A.

\sub{Proof of Theorem B.}
The spectral flow, as well as every its multiple, satisfies the first condition of the theorem.
Conversely, suppose that an additive homotopy invariant $\Phi$ 
vanishes on $\Omega_g\Ellt^0(2k_M)$ for every $k\in\N$ and $g\in U(2k_M)$. 
Then, as was stated above, $\Phi$ satisfies also the first condition of Theorem C.
It follows that there is an (invertible) constant $\lambda\in\Lambda$
such that $\Phi(\gamma)  = \Psi(\gamma,g)\cdot \lambda$ 
for every $\gamma\in\Omega_g\Ellt(2k_M)$, $k$, and $g$. 
Substituting the value of $\Psi$ given by Theorem A, $\Psi(\gamma,g) = \spf(\gamma)$,
we obtain the second condition of the theorem.

\sub{Continuity criteria.}
In the Appendix we give a general criterion of being graph continuous 
for families of closed operators in Hilbert and Banach spaces; see Proposition \ref{prop:closed}.
Then we apply this criterion to elliptic boundary value problems. Proposition \ref{prop:diff} gives a sufficient condition for continuity 
of $d$-th order differential operators with general boundary conditions. 
Its particular case concerning first order differential operators
is used in Proposition \ref{prop:AL2} and Lemma \ref{lem:lambda} in the main part of the paper.

\sub{Motivation.}
The computation of the spectral flow for paths of first order self-adjoint elliptic operators over a surface
is important for some applications in condensed matter physics.
For example, the Aharonov-Bohm effect for a single-layer graphene sheet with holes arises
if a one-parameter family of Dirac operators has non-zero spectral flow.
The varying free term of the Dirac operator corresponds to a varying magnetic field,
while the path connecting two gauge equivalent operators corresponds to the situation
where magnetic fluxes through holes change by integer numbers in the units of the flux quantum.
The spectral flow for such paths of Dirac operators was computed by the author in \cite{Pr}.
Later the results of \cite{Pr} were improved and generalized
to more general families of Dirac operators with local boundary conditions:
for even-dimensional compact manifolds by A.~Gorokhovsky and M.~Lesch in \cite{GL},
and for compact manifolds of arbitrary dimension by M.~Katsnelson and V.~Nazaikinskii in \cite{KN}.
Unfortunately, the methods of both \cite{GL} and \cite{KN} use essentially the specific nature of Dirac operators
and cannot be applied to general self-adjoint elliptic differential operators.
However, some other possible realizations of the Aharonov-Bohm effect in condensed matter physics 
are described by self-adjoint elliptic operators of non-Dirac type.
The initial motivation of the author in the present paper
was to solve the arising mathematical problem,
namely to compute the spectral flow for such a family.

\sub{Generalization to arbitrary families.}
In the next paper \cite{Pr18} we generalize Theorems A, B, and C 
to families of self-adjoint elliptic local boundary value problems 
on a compact surface parametrized by points of an arbitrary compact space $X$.
The spectral flow is then replaced by the analytical index 
and the integer-valued invariant $\Psi$ is replaced by the topological index. 
Both the analytical and the topological index take value in $K^1(X)$.
The idea of proofs remains essentially the same, 
but constant loops are replaced by \q{locally constant} families of boundary value problems,
that is, fixed boundary value problems twisted by vector bundles over $X$.

\sub{Acknowledgments.}
In the initial stages of the work reported in this paper 
I enjoyed the hospitality and excellent working conditions at the Max Planck Institute for Mathematics at Bonn. 
This work is a part of my PhD thesis at the Technion -- Israel Institute of Technology. 
I would like to use this opportunity to express my gratitude to my PhD advisor S.~Reich. 
I am also grateful to N.V.~Ivanov for his support and interest in this paper 
and to an anonymous referee for his/her suggestions.

\section{Spectral flow for unbounded operators}\label{sec:sfR}

\upskip

\sub{The space of regular operators.}
Recall that an unbounded operator $A$ on $H$ is a linear operator 
defined on a subspace $\D$ of $H$ and taking values in $H$;
the subspace $\D$ is called the domain of $A$ and is denoted by $\dom(A)$.
An unbounded operator $A$ is called closed if its graph is closed in $H\oplus H$ 
and densely defined if its domain is dense in $H$.
It is called \textit{regular} if it is closed and densely defined.

Associating with a closed operator on $H$ the orthogonal projection on its graph
defines an inclusion of the set of closed operators on $H$ into the space 
$\Proj(H\oplus H)\subset\BB(H\oplus H)$ of projections in $H\oplus H$.
Let $\Reg(H)$ be the set of regular operators on $H$ together with the topology
induced from the norm topology on $\Proj(H\oplus H)$ by this inclusion. 
This topology is usually called the \textit{graph topology}, or \textit{gap topology}.
On the subset $\B(H)\subset \Reg(H)$ it coincides with the usual norm topology 
\cite[Addendum, Theorem 1]{CL}.
So, $\B(H)$ is a subspace of $\Reg(H)$; 
it is open and dense in $\Reg(H)$ \cite[Proposition 4.1]{BLP-04}.

\sub{Fredholm operators and the spectral flow.}
Denote by $\FR(H)$ the subspace of $\Reg(H)$ consisting of regular Fredholm operators,
and by $\FR\sa(H)$ its subspace consisting of regular Fredholm self-adjoint operators.
The space $\FR\sa(H)$ is path-connected and its fundamental group is isomorphic to $\Z$ \cite{Jo}.
This isomorphism is given by the 1-cocycle on $\FR\sa(H)$ called the spectral flow.  
The definitions of the spectral flow can be found in \cite{Ph} for the case of bounded operators
and in \cite{BLP-04,L-04} for the case of unbounded operators.

The case where one or both of the endpoints of the path have zero eigenvalue 
requires some agreement on the counting procedure.
Yet if a path is a loop up to an automorphism of $H$,
the value of the spectral flow is independent of the choice of definition.
Since we consider only such paths in this paper,
we do not specify the counting agreement for the case of non-invertible endpoints: any such agreement will suffice.

\sub{Properties of the spectral flow.}
It is well known that the spectral flow has a number of nice properties:

\textbf{(S0) Zero crossing.}
In the absence of zero crossing the spectral flow vanishes:
if $\gamma$ is a continuous path in $\FR\sa(H)$ such that none of the operators $\gamma(t)$ has zero eigenvalue,
then $\sp{\gamma} = 0$.

\textbf{(S0')}
The spectral flow of a constant path vanishes.

\textbf{(S1) Homotopy invariance.}
The spectral flow along a continuous path $\gamma$ in $\FR\sa(H)$ does not change if
$\gamma$ changes continuously in the space of paths in $\FR\sa(H)$ with
fixed endpoints (the same as the endpoints of $\gamma$).

\textbf{(S2) Additivity with respect to direct sum.}
Let $H_1$, $H_2$ be separable Hilbert spaces,
and let $\gamma_i \colon [a,b] \to \FR\sa(H_i)$ be continuous paths.
Then $\sp{\gamma_1\oplus\gamma_2} = \sp{\gamma_1}+\sp{\gamma_2}$,
where $\gamma_1\oplus\gamma_2 \colon [a,b] \to \FR\sa(H_1\oplus H_2)$ denotes the pointwise direct sum.

\textbf{(S3) Path additivity.}
Let $\gamma$, $\gamma'$ be continuous paths in $\FR\sa(H)$ such that
the last point of $\gamma$ is the first point of $\gamma'$.
Then $\sp{\gamma.\gamma'} = \sp{\gamma}+\sp{\gamma'}$,
where $\gamma.\gamma'$ denotes the concatenation of $\gamma$ and $\gamma'	$.

\textbf{(S4) Conjugacy invariance.}
Let $g$ be a unitary automorphism of $H$,
and let $\gamma$ be a continuous path in $\FR\sa(H)$.
Then $\sp{\gamma} = \sp{g\gamma g\inv}$.

\sub{Paths with conjugate ends.}
In this paper we compute the spectral flow only for paths with conjugate ends (in particular, for loops),
so it is convenient to have a special designation for the space of such paths.
For a topological space $X$ we denote by $\Omega X$ the space of free loops in $X$
with the compact-open topology.
Here by a free loop we mean a continuous map from a circle $S^1$ to $X$, or, equivalently,
a continuous map $\gamma\colon [0,1] \to X$ such that $\gamma(0)=\gamma(1)$.
If $g$ is a homeomorphism of $X$, 
then we denote by $\Omega_g X$ the space of continuous paths $\gamma\colon [0,1] \to X$
such that $\gamma(1)=g\gamma(0)$ equipped with the compact-open topology.
We say that paths $\gamma, \gamma' \in\Omega_g X$ are homotopic if they can be connected by a path in $\Omega_g X$.

The group $U(H)$ of unitary automorphisms of $H$ acts on the space $\FR\sa(H)$ by conjugations:
$(A,g)\mapsto gAg\inv$.
We will write $\Omega_g\FR\sa(H)$ for $g\in U(H)$ having in mind this action.

In the proof of Theorem \ref{Thm:sf} we do not use all properties (S0-S4), 
but only the following small part of them.

\begin{prop}\label{prop:SiU}
The spectral flow has the following properties.
\begin{enumerate}\upskip
	\item[($S0_U$)] \textbf{Zero crossing.}
	Let $\gamma\in\Omega_g\FR\sa(H)$, $g\in U(H)$.
Suppose that $\gamma(t)$ has no zero eigenvalue for each $t\in[0,1]$.
Then $\spf(\gamma)=0$.
  \item[($S1_U$)] \textbf{Homotopy invariance.} 
	The spectral flow is constant on path connected components of $\Omega_g\FR\sa(H)$ for each $g\in U(H)$.
	\item[($S2_U$)] \textbf{Additivity with respect to direct sum.}
	 Let $\gamma_i\in\Omega_{g_i}\FR\sa(H_i)$, $g_i\in U(H_i)$, $i=1,2$.
Then $\spf(\gamma_1\oplus\gamma_2) = \spf(\gamma_1) + \spf(\gamma_2)$.
\end{enumerate}
\end{prop}

\proof
Properties ($S0_U$) and ($S2_U$) are just weaker versions of (S0) and (S2) respectively.
To prove ($S1_U$), we combine (S1), (S3), (S4), and (S0').
Let $\gamma_s(t)$, $s\in[0,1]$ be a homotopy between $\gamma_0$ and $\gamma_1$.
Let the paths $\beta,\beta',\beta''\colon[0,1]\to\FR\sa(H)$ be given by the formulas
$\beta(s)=\gamma_s(0)$, $\beta'(t)=\gamma_1(t)$, $\beta''(s)=\gamma_{1-s}(1)$.
Then $\gamma_0$ is homotopic to $\beta.\beta'.\beta''$ in the space of paths in $\FR\sa(H)$ 
with the same endpoints as $\gamma_0$.
Property (S1) implies $\sp{\gamma_0} = \sp{\beta.\beta'.\beta''}$,
and by (S3) the last value is equal to $\sp{\beta}+\sp{\beta'}+\sp{\beta''}$.
Property (S4) implies $\sp{\beta}=\sp{g\beta}$.
The path $g\beta$ is just the path $\beta''$ passing in the opposite direction,
so the concatenation of these two paths is homotopic to the constant path (in the class of paths with fixed endpoints).
By (S3), (S1), and (S0') we have $\sp{g\beta}+\sp{\beta''} = \sp{g\beta.\beta''} = 0$.
Taking all this together, we obtain
$\sp{\gamma_0} = \sp{\beta'} = \sp{\gamma_1}$.
\endproof

\section{Local boundary value problems}\label{sec:bound_pr}

This section is mostly devoted to standard facts about elliptic boundary value problems 
on compact manifolds,
in the context of first order operators and local boundary conditions. 
In the end of the section we give a criterion of continuity 
for families of boundary value problems,
which is a particular case of results obtained in the Appendix.

\sub{Operators.}
Let $M$ be a smooth compact connected oriented manifold with non-empty boundary $\pM$ and a fixed Riemannian metric,
and let $E$ be a smooth Hermitian complex vector bundle over $M$.
We denote by $\EY$ the restriction of $E$ to the boundary $\pM$.

Let $A$ be a first order elliptic differential operator acting on sections of $E$.
Recall that an operator $A$ is called elliptic if its (principal) symbol
$\sigma_A(\xi)$ is non-degenerate for every non-zero cotangent vector $\xi\in T^{\ast}M$.
Throughout the main part of the paper (except for the Appendix) all differential operators are supposed 
to have smooth ($C^{\infty}$) coefficients.

Green's formula for $A$ can be written as 
\begin{equation}\label{eq:Green}
  \bra{Au,v}_{L^2(E)} - \bra{u,A^t \, v}_{L^2(E)} = \bra{i\sigma_A(n)\restr{u}{\pM},\restr{v}{\pM}}_{L^2(\EY)}
  \; \text{ for } u,v\in C^{\infty}(E), 
\end{equation}
where $A^t$ denotes the differential operator formally adjoint to $A$.
If $A$ is formally self-adjoint, then it is symmetric on the domain $C^{\infty}_0(E)$,
that is $\bra{Au,v}_{L^2(E)} = \bra{u,Av}_{L^2(E)}$
for any smooth sections $u$, $v$ of $E$ with compact supports in $M\setminus\pM$.

\sub{Local boundary conditions.}
The differential operator $A$ with the domain $C_0^{\infty}(E)$
is an unbounded operator on the Hilbert space $L^2(E)$ of $L^2$-sections of $E$.
This operator can be extended to a closed operator on $L^2(E)$ in various ways, 
by imposing appropriate boundary conditions.
We will consider only local boundary conditions that are defined by smooth subbundles of $\EY$.
For such a subbundle $L$, the corresponding unbounded operator $A_L$ on $L^2(E)$ has the domain
\begin{equation}\label{eq:dom}
  \dom\br{A_L} = \set{ u\in H^1(E)\colon \restr{u}{\pM} \mbox{ is a section of } L },	
\end{equation}
where $H^1(E)$ denotes the first order Sobolev space
(the space of sections of $E$ which are in $L^2$ together with all their first derivatives).
We will often identify a pair $(A,L)$ with the operator $A_L$.

To give a precise meaning to the notation in the right-hand side of \eqref{eq:dom},
recall that the restriction map $\Cinf(E)\to \Cinf(\EY)$ 
taking a section $u$ to $\restr{u}{\pM}$
extends continuously to the trace map $\tau\colon H^1(E)\to H^{1/2}(\EY)$.
The smooth embedding $L\hookto\EY$ defines the natural inclusion $H^{1/2}(L)\hookto H^{1/2}(\EY)$.
By the condition \q{$\restr{u}{\pM}$ is a section of $L$} in \eqref{eq:dom}
we mean that the trace $\tau(u)$ lies in the image of this inclusion.

\sub{Decomposition of $E$.}
We will use the following properties of elliptic symbols.

\begin{prop}\label{prop:Qn}
Let $\sigma\in\Hom(T^*M,\End(E))$ be a symbol of first order elliptic operator.
Let $\Pi$ be an oriented two-dimensional plane in the cotangent bundle $T^*_x M$, $x\in M$.
Then for any positive oriented frame $(e_1, e_2)$ in $\Pi$
the operator 
$$Q = \sigma(e_1)^{-1} \sigma(e_2)\in \End(E_x)$$
has no eigenvalues on the real axis.
It defines the direct sum decomposition $E_x = \Ep \oplus \Em$ (not necessarily orthogonal), 
where $\Ep$ and $\Em$ are the generalized eigenspaces of $Q$
corresponding to the eigenvalues with positive and negative imaginary part respectively.
This decomposition depends only on $\Pi$ and is independent of the choice of a frame $(e_1, e_2)$.

If additionally $\sigma$ is self-adjoint, 
then the ranks of $\Ep$ and $\Em$ are equal (so the rank of $E$ is even), 
and for every non-zero $\xi\in\Pi$ 
the symbol $\sigma(\xi)$ takes $\Ep$ and $\Em$ to their orthogonal complements in $E_x$.
\end{prop}

\proof
1. Since $\sigma$ is elliptic, the operator
$Q-t = \sigma(e_1)\inv \sigma(e_2-te_1)$ is invertible for any $t\in\R$.
Hence $Q$ has no eigenvalues on the real axis and $E_x = \Ep\oplus\Em$.

If we change $(e_1,e_2)$ to $(e_1, e_2+te_1)$, $t\in\R$, then $Q$ is changed to $Q+t\Id$.
If we change $(e_1,e_2)$ to $(e_1+te_2, e_2)$ then $Q$ is changed to $(Q\inv+t\Id)\inv$.
In both cases $\Ep$ and $\Em$ do not change.
Therefore, they do not change at any change of the frame $(e_1,e_2)$ preserving orientation,
and thus depend only on $\Pi$.

2. Suppose now that $\sigma$ is self-adjoint, that is $\sigma(\xi)$ is self-adjoint for every $\xi\in T^* M$.
Let $\xi\in\Pi$ be a non-zero vector.
Choose a positive oriented frame $(e_1, e_2)$ in $\Pi$ such that $e_1=\xi$.	
Denote $\sigma_i = \sigma(e_i)$, $V_{\lambda,\,k} = \Ker(Q-\lambda)^k$, and $V_{\lambda} = V_{\lambda,\,\dim E}$.
We prove by induction that $\sigma_1 V_{\lambda}$ is orthogonal to $V_{\mu}$
for any $\lambda, \mu\in\CC$ with $\lambda\neq\overline{\mu}$. 
Indeed, $\sigma_1 V_{\lambda,\,0} = 0$ is orthogonal to $V_{\mu,\,0} = 0$.
Suppose that $\sigma_{1}V_{\lambda,\,l}$ is orthogonal to $V_{\mu,\,m}$ for all $l,m\geq 0$, $l+m<k$.
Then for $l+m=k$, $u\in V_{\lambda,\,l}$, $v\in V_{\mu,\,m}$ we have
\begin{multline*}
(\lambda-\overline{\mu}) \bra{\sigma_1 u,v} =
   \bra{\sigma_1\lambda u, v} - \bra{\sigma_1 u, \mu v} + \bra{u, \sigma_2 v} - \bra{\sigma_2 u,v} =
	\\
	=  \bra{\sigma_1\lambda u, v} - \bra{\sigma_1 u, \mu v} + \bra{\sigma_1 u, Qv} - \bra{\sigma_1 Qu,v} =
	\bra{\sigma_1 u, (Q-\mu)v} - \bra{\sigma_1(Q-\lambda)u, v} = 0
\end{multline*}
by induction assumption, since $(Q-\mu)v\in V_{\mu,\,m-1}$ and $(Q-\lambda)u\in V_{\lambda,\,l-1}$.
Thus $\sigma_{1}V_{\lambda}$ is orthogonal to $V_{\mu}$ if $\lambda\neq\overline{\mu}$.

The subspace $\Ep$ is spanned by $\bigcup V_{\lambda}$ with $\lambda$ running all eigenvalues of $Q$ with positive imaginary parts.
For every pair $\lambda$, $\mu$ of such eigenvalues (not necessarily distinct) we have
$\lambda\neq\overline{\mu}$, so
$\sigma_{1}\Ep$ is orthogonal to $\Ep$.
Similarly, $\sigma_{1}\Em$ is orthogonal to $\Em$.
We have
$$2\dim\Ep = \dim\Ep + \dim(\sigma_{1}\Ep) \leq \dim\Ep + \dim(\Ep)^{\bot} = \dim E_x$$
and, similarly, $2\dim\Em \leq \dim E_x$.
On the other hand, $\dim\Ep+\dim\Em = \dim E_x$.
Therefore, $\dim\Ep = \dim\Em = \dim E_x/2$.
\endproof

\sub{Elliptic boundary conditions.}
Let $A$ be a first order elliptic operator acting on sections of $E$.
The inverse image of a non-zero cotangent vector $\xi\in T^*_x\pM$ 
under the restriction map $T^*_x M\to T^*_x\pM$ 
is an affine line in $T^*_x M$ parallel to the outward conormal $n_x$. 
Denote by $\Pi_{\xi}$ the two-dimensional vector subspace of $T^*_x M$ spanned by this line.
Applying Proposition \ref{prop:Qn} to the plain $\Pi_{\xi}$, 
we obtain the decomposition of $E_x$ into the direct sum $\Ep(\xi)\oplus\Em(\xi)$. 

A local boundary condition $L$ is called elliptic for $A$ if
\begin{equation}\label{eq:Lell-gen}
	L_x \cap \Ep(\xi) = 0 \;\text{ and }\; L_x+\Ep(\xi) = E_x 
	\;\text{ for every non-zero } \xi\in T_x^*\pM, \; x\in\pM.
\end{equation}
If $L$ is elliptic for $A$, then the adjoint to $A_L$ is $A^t_{N}$, 
where $A^t$ is the differential operator formally adjoint to $A$ 
and $N = (\sigma(n)L)^{\bot}$.

\sub{Self-adjoint elliptic boundary conditions.}
Suppose now that an elliptic operator $A$ is formally self-adjoint.
Then the conormal symbol $\sigma(n)$ of $A$ defines a symplectic structure on fibers of $\EY$
given by the symplectic 2-form
\[
\omega_x\colon E_x\otimes E_x\to\CC, \;
\omega_x(u,v) = \bra{i\sigma(n_x) u, \;v} \; \text{ for } u,v\in E_x, 
\]
where $n_x$ is the outward conormal to $\pM$ at $x\in \pM$.
By Proposition \ref{prop:Qn} both $\Ep(\xi)$ and $\Em(\xi)$ 
are Lagrangian subspaces with respect to this symplectic structure.

The differential operator $A$ with the domain $C^{\infty}_0(E)$
is a symmetric unbounded operator on the Hilbert space $L^2(E)$ of $L^2$-sections of $E$.
This operator can be extended to a regular self-adjoint operator
on $L^2(E)$ by imposing appropriate boundary conditions.
For $L$ satisfying ellipticity condition \eqref{eq:Lell-gen},
the operator $A_L$ is self-adjoint if and only if $L$ is a Lagrangian subbundle of $\EY$. 

For a Lagrangian subbundle $L$ condition \eqref{eq:Lell-gen} can be written in simpler form:
\begin{equation}\label{eq:Lell-gen2}
	L_x \cap \Ep(\xi) = 0 \;\text{ for every non-zero } \xi\in T_x^*\pM, x\in\pM.
\end{equation}
Indeed, rank of both $L_x$ and $\Ep(\xi)$ is half of rank $E_x$. 
Therefore, $L_x \cap \Ep(\xi) = 0$ if and only if $L_x+\Ep(\xi) = E_x$.

Finally, we obtain the following description of self-adjoint elliptic local boundary value problems.

\begin{prop}\label{prop:AL1}
Let $A$ be a first order formally self-adjoint elliptic differential operator acting on sections of $E$.
Let $L$ be a smooth Lagrangian subbundle of $\EY$ satisfying condition \eqref{eq:Lell-gen2}. 
Then $A_L$ is a regular Fredholm self-adjoint operator on $L^2(E)$.
Moreover, $A_L$ has compact resolvents, that is $(A_L+i)\inv$ is a compact operator on $L^2(E)$.
\end{prop}

\proof
Denote by $\D$ the domain of $A_L$ given by formula \eqref{eq:dom}.
It is dense in $L^2(E)$ and closed in $H^1(E)$. 
Equip $\D$ with the topology induced from $H^1(E)$.

Let $\tau\colon H^1(E)\to H^{1/2}(\EY)$ be the trace map, 
and let $P$ be the bundle endomorphism projecting $\EY$ on $L^{\bot}$ along $L$.
Condition \eqref{eq:Lell-gen} means that $P\colon\Ep(\xi)\to L_x^{\bot}$ is bijective
for every non-zero $\xi\in T_x^*\pM$.
It follows by \cite[Theorem 20.1.2]{Hor} that the operator 
$A\oplus P\tau \colon H^1(E) \to L^2(E)\oplus H^{1/2}(L^{\bot})$ is Fredholm. 
Its restriction to the kernel of $P\tau$ is also Fredholm.
But this restriction coincides with $A_L$ considered as a bounded operator from $\D$ to $L^2(E)$. 
Hence $A_L$ is Fredholm.

In particular, $V=\im(A_L)$ is a closed subspace of $L^2(E)$.
Let $U$ be the orthogonal complement of the kernel of $A_L$ in $L^2(E)$.
The restriction $\bar{A}$ of $A$ to $U$ is injective with the image $V$.
Therefore, the inverse operator $\bar{A}\inv\colon V\to U$ is bounded and its graph is closed in $V\times U$.
Equivalently, the graph of $\bar{A}$ is closed in $U\times V$, which is a closed subspace of $L^2(E)^2$.
The graph of $A_L$ is the orthogonal sum of $\ker(A_L)\times\set{0}$ with the graph of $\bar{A}$ 
and therefore is closed in $L^2(E)^2$.
In other words, the operator $A_L$ is closed.

Green's formula \eqref{eq:Green}
implies that $A_L$ is symmetric.
Let $(u,v)\in L^2(E)^2$ be an arbitrary point of the graph of the adjoint operator.
This means that for each $w\in\dom(A_L)$ we have $\bra{u,Aw}=\bra{v,w}$. 
By \cite[Theorem 1]{Lax}, $(u,v)$ lies in the closure of the graph of $A_L$.
(The statement of this theorem of Lax and Phillips concerns only 
smooth domains in Euclidean spaces and trivial vector bundles. 
But its proof is local, so it works for general case without change.)
Since $A_L$ is closed, $u\in\dom(A_L)$.
Therefore, $A_L$ is self-adjoint.

The operator $A_L$ is bounded as an operator from a Hilbert space $\D$ to $L^2(E)$.
Since $A_L$ is a closed self-adjoint operator on $L^2(E)$,
the bounded operator $A_L+i\colon\D\to L^2(E)$ is bijective and
the inverse $(A_L+i)\inv$ is a bounded operator from $L^2(E)$ to $\D$ \cite[Theorem V.3.16]{Kato}.
Composing it with the compact embedding $\D\subset H^1(E)\hookto L^2(E)$, 
we see that $(A_L+i)\inv$ is compact as an operator on $L^2(E)$.
This completes the proof of the proposition.
\endproof

\sub{The space of boundary value problems.}
Denote by $\Ellt(E)$ the set of all pairs $(A,L)$ satisfying conditions of Proposition \ref{prop:AL1}. 
The following result is a particular case of Proposition \ref{prop:AL2-app} from the Appendix.

\begin{prop}\label{prop:AL2}
For the set $\Ellt(E)$ equipped with the $C^0$-topology on coefficients of operators 
and the $C^1$-topology on boundary conditions,
the natural inclusion $\Ellt(E) \hookto \FR\sa\br{L^2(E)}$, $(A,L) \mapsto A_L$ is continuous.
\end{prop}

\upskip
Equivalently, the $C^0$-topology on coefficients of operators 
can be described as the topology induced by the inclusion $\Ell(E) \hookto \BB(H^1(E), L^2(E))$.

\section{Boundary value problems on a surface}\label{sec:bound_pr2}

From now on (except for the Appendix) we will consider only the case of \textbf{\textit{dimension two}}, 
that is $M$ will be a smooth compact connected oriented \textbf{\textit{surface}} 
with non-empty boundary $\pM$ and a fixed Riemannian metric.

Let $E$ be a smooth Hermitian complex vector bundle over $M$.
Denote by $\Ell(E)$ the set of first order formally self-adjoint elliptic differential operators 
with smooth coefficients acting on sections of $E$.

\sub{Decomposition of a bundle.}
Since $M$ is now two-dimensional, 
Proposition \ref{prop:Qn} allows to define the \textit{global} decomposition of $E$.

\begin{prop}\label{prop:Q2}
Let $A\in\Ell(E)$. 
Then the symbol $\sigma$ of $A$ defines the decomposition of $E$ into the direct sum (not necessarily orthogonal) 
of two smooth subbundles $\Ep=\Ep(\sigma)$ and $\Em=\Em(\sigma)$ such that the following conditions hold:
\begin{enumerate}\upskip
	\item $\Ep_x$ and $\Em_x$ are the generalized eigenspaces of $Q_x = \sigma(e_1)^{-1} \sigma(e_2)$ as in Proposition \ref{prop:Qn},
where $(e_1, e_2)$ is an arbitrary positive oriented frame in $T^{*}_x M$.
	\item Ranks of $\Ep$ and $\Em$ are equal, so the rank of $E$ is even.
	\item For every non-zero $\xi\in T^{*}_x M$ the symbol $\sigma(\xi)$ takes $\Ep_x$ and $\Em_x$ to their orthogonal complements in $E_x$.
\end{enumerate}\upskip
\end{prop}

\proof
The main part of the statement follows from Proposition \ref{prop:Qn}.
It remains to show that $\Ep_x$ and $\Em_x$ are fibers of smooth vector bundles $\Ep$ and $\Em$.
Choosing a local smooth frame $(e_1, e_2)$ in $T^{*} M$,
we see that $\Ep_x$ and $\Em_x$ smoothly depend on $Q_x$, which in turn smoothly depends on $x$.
\endproof

\sub{Self-adjoint elliptic boundary conditions.}
Denote by $\EpY$, resp. $\EmY$ the restriction of $\Ep$, resp. $\Em$ to $\pM$.
As before, the conormal symbol $\sigma(n)$ defines the symplectic structure on the fibers of $\EY$,
and $\EpY$, $\EmY$ are transversal Lagrangian subbundles of $\EY$.

The orientation on $M$ induces the orientation on $\pM$.
Fibers of $\EY^{\pm}$ can be written as $\Ep_x = \Ep(\xi)$ and $\Em_x = \Em(\xi)$,
where $\xi$ is a positive vector in the oriented one-dimensional space $T^{*}_x \pM$.
The identity $\Ep(\xi) = \Em(-\xi)$ allows to write ellipticity condition \eqref{eq:Lell-gen} in simpler form:
\begin{equation}\label{eq:L1dim2}
	L \cap \EpY = L \cap \EmY = 0 \;\text{ and }\; L+\EpY = L+\EmY = \EY.
\end{equation}
If $L$ is Lagrangian, then condition \eqref{eq:L1dim2} can be simplified even further,
cf. \eqref{eq:Lell-gen2}:
\begin{equation*}
	L \cap \EpY = L \cap \EmY = 0.
\end{equation*}
As before, we denote by $\Ellt(E)$ the set of all pairs $(A,L)$ such that $A\in\Ell(E)$
and $L$ is a smooth Lagrangian subbundle of $\EY$ satisfying condition \eqref{eq:L1dim2}. 
Proposition \ref{prop:AL1} then takes the following form. 

\begin{prop}\label{prop:AL1dim2}
For every $(A,L)\in\Ellt(E)$ the unbounded operator $A_L$ 
is a regular self-adjoint operator on $L^2(E)$ with compact resolvents.
\end{prop}

\upskip

\sub{The correspondence between boundary conditions and automorphisms of $\EmY$.}
For every elliptic (not necessarily self-adjoint) symbol $\sigma$ 
there is a one-to-one correspondence between 
subbundles $L$ of $\EY$ satisfying condition \eqref{eq:L1dim2} and
bundle isomorphisms $R\colon\EmY\to\EpY$.
Namely, $L$ is the graph of $R$ in $\EmY\oplus\EpY = \EY$.
Equivalently, $-R$ is the projection of $\EmY$ onto $\EpY$ along $L$.

If additionally $\sigma$ is self-adjoint, then one can move further and construct 
a one-to-one correspondence between \textit{Lagrangian} subbundles $L$ of $\EY$ satisfying \eqref{eq:L1dim2} 
and \textit{self-adjoint} bundle automorphisms $T$ of $\EmY$.

Let us first describe this correspondence in the case of \textit{mutually orthogonal} $\EpY$ and $\EmY$ 
(this holds, in particular, for Dirac type operators).
Composing $R$ with $i\sigma(n)\colon\EpY\to(\EpY)^{\bot}=\EmY$, we obtain the bundle automorphism $T$ of $\EmY$.
Conversely, with every bundle automorphism $T$ of $\EmY$ we associate the subbundle $L$ of $\EY$ 
given by the formula 
\begin{equation}\label{eq:ortL}
		L=\set{ u^+\oplus u^-\in \EpY\oplus\EmY = \EY \colon \; i\sigma(n)u^+=Tu^- }.
\end{equation}
As Proposition \ref{prop:T} below shows, 
$T$ is self-adjoint if and only if $L$ is Lagrangian, 
so we obtain a bijection between the set of all self-adjoint elliptic local boundary conditions for $A$
and the set of all self-adjoint bundle automorphisms of $\EmY$.

In general case, where $\EpY$ and $\EmY$ can be non-orthogonal, 
this construction should be slightly modified.
The composition $\tilde{T} = i\sigma(n)R$ acts from $\EmY$ to $(\EpY)^{\bot}$, 
which now does not coincide with $\EmY$.
In order to correct this, we compose $\tilde{T}$ with the orthogonal projection 
$\Pmort$ of $\EY$ onto $\EmY$.
Since $\EY = (\EmY)^{\bot}\oplus(\EpY)^{\bot}$, 
the restriction of $\Pmort$ to $(\EpY)^{\bot}$ is an isomorphism $(\EpY)^{\bot}\to\EmY$.
Finally, we define the bundle automorphism $T = \Pmort\circ\tilde{T}$ of $\EmY$,
so that the following diagram becomes commutative.
\begin{equation}\label{diag:T}							
\begin{tikzcd}
L \arrow{d}[swap]{P^-} \arrow{r}{P^+} & \EpY \arrow{r}{i\sigma(n)} 
& (\EpY)^{\bot} \arrow[swap]{d}{\Pmort} \\
\EmY \arrow[dashed]{ru}{R} \arrow[swap, dashed]{rru}{\tilde{T}} \arrow[swap,dashed]{rr}{T} 
&& \EmY \arrow[swap, shift right=1ex]{u}{(P^-)^{\ast}} 
\end{tikzcd}
\end{equation}

\begin{prop}\label{prop:T}
Let $A\in\Ell(E)$.
Denote by $P^+$ the projection of $\EY$ onto $\EpY$ along $\EmY$
and by $P^- = 1-P^+$ the projection of $\EY$ onto $\EmY$ along $\EpY$.
Then the following hold.
\begin{enumerate}\upskip
	\item There is a one-to-one correspondence between 
	      smooth subbundles $L$ of $\EY$ satisfying condition \eqref{eq:L1dim2} and
				smooth bundle automorphisms $T$ of $\EmY$.
				This correspondence is given by the formula
	      \begin{equation}\label{eq:PT}
	        L=\Ker P_T \; \mbox{ with } \; P_T = P^+\br{ 1+i\sigma(n)\inv TP^- },
	      \end{equation}
	      where $P_T$ is the projection of $\EY$ onto $\EpY$ along $L$.
	\item For $L$ and $T$ as above, $L$ is Lagrangian if and only if $T$ is self-adjoint.
\end{enumerate}
\upskip
If $\EpY$ and $\EmY$ are mutually orthogonal, then \eqref{eq:PT} is equivalent to \eqref{eq:ortL}.
\end{prop}

\noindent
In the rest of the paper
we will sometimes write an element of $\Ellt(E)$ as $(A,T)$ instead of $(A,L)$.

\medskip

\proof
The adjoint $(P^-)^{\ast}$ projects $\EY$ onto $(\EpY)^{\bot}$ along $(\EmY)^{\bot}$,
so its restriction to $\EmY$ is the inverse of $\Pmort\colon(\EpY)^{\bot}\to\EmY$.
All three solid arrows at the right half of Diagram \eqref{diag:T} are smooth bundle isomorphisms.

By Proposition \ref{prop:Q2} the conormal symbol $\sigma(n)$ 
takes $\EmY$ and $\EpY$ to their orthogonal complements.
So $(P^-)^{\ast} = \sigma(n)P^+\sigma(n)\inv$ and
$\sigma(n)\inv(P^-)^{\ast} = P^+\sigma(n)\inv$.
Therefore, \eqref{eq:PT} can be equivalently written as 
\begin{equation}\label{eq:PT2}
	P_T = P^+ + i\sigma(n)\inv(P^-)^{\ast}TP^-.
\end{equation}

\noindent
1. Let $L$ be a smooth subbundle of $\EY$ satisfying \eqref{eq:L1dim2}.
Then both solid arrows at the left half of Diagram \eqref{diag:T} are smooth bundle isomorphisms.	
There is a smooth automorphism $T$ of $\EmY$ making this diagram commutative,
and such an automorphism is unique.
Substituting $R = (i\sigma(n))\inv(P^-)^{\ast}T$ to $L = \Ker(P^+-RP^-)$,
we obtain $L = \Ker(P^+ + i\sigma(n)\inv(P^-)^{\ast}TP^-) = \Ker P_T$.

Conversely, let $T$ be a smooth automorphism of $\EmY$.
The image of $P_T$ is contained in $\EpY$, while the restriction of $P_T$ to $\EpY$ is the identity.
It follows that $P_T^2=P_T$, that is
$P_T$ is the projection of $\EY$ onto $\EpY$ along $L=\Ker P_T$.
This implies $L\cap\EpY = 0$ and $L+\EpY =\EY$.
The projection $P_T$ smoothly depends on $x\in\pM$ and has constant rank, 
so $L$ is a smooth subbundle of $\EY$ with $\rank L = \rank\EY - \rank\EpY = \rank\EmY$.
If $u\in L\cap\EmY$, then $P^+ u = 0$ and
$Tu = \Pmort i\sigma(n)P^+ u = 0$.
Since $T$ is invertible, $L\cap\EmY = 0$.
This completes the proof of clause 1.

\noindent
2. Let $L$, $T$ be as in clause 1 and $u_1,u_2\in L$.
For $u_j^- = P^-u_j$ and $u_j^+ = P^+u_j$ we have
\begin{equation}\label{eq:Tuv1}
  \bra{Tu_1^-,u_2^-} = \bra{\tilde{T}u_1^-,u_2^-} =	\bra{i\sigma(n)u_1^+, u_2^-} = \bra{i\sigma(n)u_1^+, u_2},
\end{equation}
using the orthogonality of $\tilde{T}u_1^- - Tu_1^- = (1-\Pmort)\tilde{T}u_1^- \in(\EmY)^{\bot}$ to $u_2^-\in\EmY$
and the orthogonality of $i\sigma(n)u_1^+\in(\EpY)^{\bot}$ to $u_2-u_2^- = u_2^+\in\EpY$.
Similarly,
\begin{equation}\label{eq:Tuv2}
  \bra{u_1^-,Tu_2^-} = \bra{u_1^-, i\sigma(n)u_2^+} = - \bra{i\sigma(n)u_1^-, u_2^+} = - \bra{i\sigma(n)u_1^-, u_2}.
\end{equation}
Subtracting \eqref{eq:Tuv2} from \eqref{eq:Tuv1}, we obtain
\begin{equation}\label{eq:Tuv}
	\bra{i\sigma(n)u_1, u_2} = \bra{TP^-u_1, P^-u_2} - \bra{P^-u_1, TP^-u_2} \;\; \text{for all} \; u_1,u_2\in L.
\end{equation}
If $L$ is Lagrangian, then \eqref{eq:Tuv} implies self-adjointness of $T$, 
since $P^-\colon L\to\EmY$ is surjective.
Conversely, if $T$ is self-adjoint, then \eqref{eq:Tuv} implies $i\sigma(n)L \subset L^{\bot}$;
taking into account that $\rank L = \rank\EY/2$, we see that $L$ is Lagrangian.

\noindent
3. If $\EpY$ and $\EmY$ are mutually orthogonal, 
then $(P^-)^{\ast}\colon\EmY\to(\EpY)^{\bot}$ is the identity,
and \eqref{eq:PT2} takes the form \eqref{eq:ortL}.
This completes the proof of the proposition.
\endproof

\sub{The subbundle $F(A,L)$.}
With every $(A,L)\in\Ellt(E)$ we associate the smooth subbundle $F(A,L)$ of $\EmY$ as follows.
Let $T$ be the self-adjoint automorphism of $\EmY$ given by formula \eqref{eq:PT}.
We define $F_x$ as the invariant subspace of $T_x$
spanned by the generalized eigenspaces of $T_x$ corresponding to negative eigenvalues.
Subspaces $F_x$ of $\Em_x$ smoothly depend on $x\in\pM$
and therefore are fibers of the smooth subbundle $F = F(A,L)$ of $\EmY$.

Being a subbundle of $\EmY$, $F(A,L)$ is also a smooth subbundle of $\EY$.
Sometimes it will be more convenient for us to consider $F(A,L)$ as a subbundle of $\EY$.

\section{The space of boundary value problems on a surface}\label{sec:Ellt}

\upskip

\sub{Topology on $\Ellt(E)$.}
In section \ref{sec:bound_pr} we used the $C^0$-topology on coefficients of operators.
We will compute the spectral flow for the paths in $\Ellt(E)$ which are continuous in a slightly stronger topology,
namely the $C^{1}$-topology on symbols and the $C^{0}$-topology on free terms of operators. 
Let us describe it more precisely.

For a smooth complex vector bundle $V$ over a smooth manifold $N$, 
we denote by $\Gr(V)$ the smooth bundle over $N$
whose fiber over $x\in N$ is the complex Grassmanian $\Gr(V_x)$.
In the same manner we define the smooth bundle $\End(V)$ of fiber endomorphisms. 
We identify sections of $\Gr(V)$ with subbundles of $V$ 
and sections of $\End(V)$ with bundle endomorphisms of $V$.

Let $r=(r_1,r_0)$ be a couple of integers, $r_1\geq r_0\geq 0$.
Denote by $\Ell^r(E)$ the set $\Ell(E)$ equipped with the $C^{r_1}$-topology on symbols
and the $C^{r_0}$-topology on free terms of operators.

To be more precise, notice that the tangent bundle $TM$ is trivial since $M$ is a surface with non-empty boundary.
Thus we can choose smooth global sections $e_1$, $e_2$ of $TM$ such that
$e_1(x)$, $e_2(x)$ are linear independent for any $x\in M$.
Choose a smooth unitary connection $\nabla$ on $E$.
Each $A\in\Ell(E)$ can be written uniquely as 				
$A = \sigma_1\nabla_1 + \sigma_2\nabla_2 + a$,	
where the symbol components $\sigma_i = \sigma_A(e_i)$ are self-adjoint bundle automorphisms of $E$,
$\nabla_i = \nabla_{e_i}$, and the free term $a$ is a bundle endomorphism.
Therefore the choice of $(e_1,e_2,\nabla)$ defines the inclusion
$$\Ell(E) \hookto \Cinf\br{\End(E)}^2 \times \Cinf\br{\End(E)},
  \quad \sigma_1\nabla_1 + \sigma_2\nabla_2 + a \mapsto \br{\sigma_1, \sigma_2, a},$$
where $\Cinf\br{\End(E)}$ denotes the space of smooth sections of $\End(E)$.
We equip $\Ell(E)$ with the topology induced by the inclusion 
$$\Ell(E) \hookto C^{r_1}\br{\End(E)}^2 \times C^{r_0}\br{\End(E)}$$
and denote the resulting space by $\Ell^r(E)$.
Equip $\Ellt(E)$ with the topology induced by the inclusion
$\Ellt(E) \hookto \Ell^r(E)\times C^{1}(\Gr(\EY))$,
with the product topology on the last space, and denote the resulting space by $\Ellt^r(E)$.
Thus defined topologies on $\Ell^r(E)$, $\Ellt^r(E)$ are independent of the choice 
of a frame $(e_1,e_2)$ and connection $\nabla$.

By Proposition \ref{prop:AL1dim2} 
the natural inclusion $\Ellt^{\,(0,0)}(E) \hookto  \FR\sa\br{L^2(E)}$ is continuous.
Since the $(r_1,r_0)$-topology on $\Ellt(E)$ is stronger than the $(0,0)$-topology, 
the inclusion $\Ellt^{r}(E) \hookto  \FR\sa\br{L^2(E)}$ is continuous for every couple $r$ of non-negative integers.

\textbf{Convention.}
\textit{From now on we will use the $(1,0)$-topology on $\Ellt(E)$,
that is, the $C^{1}$-topology on symbols and the $C^{0}$-topology on free terms of operators.
For brevity we will omit the superscript, so further $\Ellt(E)$ will always mean $\Ellt^{(1,0)}(E)$.
}

The following is an immediate corollary of Proposition \ref{prop:AL1dim2}.

\begin{prop}\label{prop:AL10}
The natural inclusion $\Ellt(E) \hookto  \FR\sa\br{L^2(E)}$, $(A,L) \mapsto A_L$ is continuous.
\end{prop}

\begin{rem}
We choose to use the stronger $(1,0)$-topology on $\Ellt(E)$ instead of the $(0,0)$-topology to simplify the proofs.
Probably, all theorems in the paper remain valid for $(0,0)$-topology on $\Ellt(E)$ as well, 
but the author did not check this.
It can be easily seen that Theorems B and C are valid (and their proofs remains the same)
for $(r_1,r_0)$-topology on $\Ellt(E)$ with $r_1-1\geq r_0\geq 0$.
\end{rem}

\upskip

\sub{Continuity of the decomposition.}
We prove here a technical result that will be used further in the paper.

Denote by $\Sigma(E)$ the set of all smooth bundle morphisms $\sigma\colon T^{\ast}M\to\End(E)$
such that $\sigma$ is a symbol of a formally self-adjoint elliptic operator.
Equip $\Sigma(E)$ with the topology induced by the inclusion
$\Sigma(E)\hookto C^1(TM\otimes\End(E))$.
Then the natural projection $\Ell(E)\to\Sigma(E)$ is continuous,
as well as the map $\Sigma(E)\to C^1(\End(\EY))$ taking $\sigma$ to $\sigma(n)$.

For a smooth fiber bundle $V$ over a smooth compact manifold $N$,
we denote by $C^{\infty,s}(V)$ the space of smooth sections of $V$ with the $C^{s}$-topology,
that is, the topology induced by the embedding $C^{\infty}(V) \hookto C^{s}(V)$.

Let $e_1$, $e_2$ be global sections of $T^{\ast}M$
such that $(e_1(x), e_2(x))$ is a positive oriented orthonormal basis of $T_x^{\ast}M$ for any $x\in M$.

\begin{prop}\label{lem:Empcont}
The following maps are continuous:
\begin{enumerate}\upskip
	\item $Q \colon \Sigma(E)\to C^{\infty,1}(\End(E))$, $\sigma\mapsto Q = \sigma(e_1)^{-1} \sigma(e_2)$;
	\item $\Ep,\Em \colon \Sigma(E) \to C^{\infty,1}(\Gr(E))$.
\end{enumerate}
\end{prop}

\proof
1. The maps from $\Sigma(E)$ to $C^{\infty,1}(\End(E))$ taking $\sigma$ to $\sigma(e_i)$, $i=1,2$, are continuous, 
so $Q$ is also continuous.

2. The invariant subspace $\Em_x$ of $Q_x$ spanned by the generalized eigenspaces of $Q_x$
corresponding to eigenvalues with negative imaginary part
is an analytic function of $Q_x$ and hence an analytic function of $\sigma_x$.
Therefore, for smooth $\sigma$, $\Em(\sigma)$ is a smooth subbundle of $E$,
and the map $\Em \colon \Sigma(E) \to C^{\infty,1}(\Gr(E))$ is continuous.
The same is true for $\Ep \colon \Sigma(E) \to C^{\infty,1}(\Gr(E))$.

\sub{Correspondence between $L$ and $T$.}
Proposition \ref{prop:T} defines a one-to-one correspondence between $L$ and $T$.
We will use it to construct homotopies in $\Omega_g\Ellt(E)$.
To do this, we need to show that the map $(A,L)\mapsto(A,T)$ is a homeomorphism.

Denote by $\Ellt'(E)$ the set of all pairs $(A,T)$ such that $A\in\Ell(E)$
and $T$ is a smooth bundle automorphism of $\EmY(A)$.
We equip $\Ellt'(E)$ with the topology induced by the inclusion 
\begin{equation}\label{eq:Ell'}
	\Ellt'(E) \hookto \Ell(E)\times C^1(\End(\EY)), \quad (A,T) \mapsto (A, T\oplus\Id_{(\EmY)^{\bot}}),
\end{equation}
where $(\EmY)^{\bot}$ is the orthogonal complement of $\EmY(A)$ in $\EY$.
We introduce the auxiliary self-adjoint automorphism 
\begin{equation}\label{eq:T'}
	T' = T\oplus\Id_{(\EmY)^{\bot}}
\end{equation}
by technical reasons:
$T$ acts on the bundle $\EmY(A)$ which depends on $A$,
while $T'$ acts on the fixed bundle $\EY$.

\begin{prop}\label{lem:T'}
	The map $\Ellt(E) \to \Ellt'(E)$ taking $(A,L)$ to $(A,T)$ is a homeomorphism.
	The map $F \colon \Ellt(E) \to C^{\infty,1}(\Gr(\EY))$ is continuous.
\end{prop}

\proof
Denote by $\Gr^{(2)}(E)$ the smooth subbundle of $\Gr(E)\times_M\Gr(E)$
whose fiber over $x\in M$ consists of pairs $(V_x,W_x)$ of subspaces of $E_x$ such that 
$V_x\cap W_x = 0$ and $V_x+W_x = E_x$.
For a smooth section $(V,W)$ of $\Gr^{(2)}(E)$ the projection $P_{V,W}$ of $E$ on $V$ along $W$ 
is a smooth section of $\End(E)$.
The map $C^{\infty,1}(\Gr^{(2)}(E)) \to C^{\infty,1}(\End(E))$, $(V,W)\to P_{V,W}$ is continuous.
The same is true if we replace $M$ by $\pM$ and $E$ by $\EY$.
Therefore, the composition
\begin{align*}
\Sigma(E) \to C^{\infty,1}(\Gr^{(2)}(E)) \to C^{\infty,1}(\Gr^{(2)}(\EY)) \to C^{\infty,1}(\End(\EY)),
\end{align*}
$\sigma \mapsto (E^+,E^-) \mapsto (\EY^+,\EY^-) \mapsto  P_{\EY^+,\,\EY^-} = P^+(\sigma)$,
is continuous.
Similarly, the map $P^-\colon \Sigma(E) \to C^{\infty,1}(\End(\EY))$ is continuous.

Let $T'$ be defined by formula \eqref{eq:T'}.
Since $TP^- = T'P^-$, identity \eqref{eq:PT2} can be equivalently written as 
$P_T = P^+ + i\sigma(n)\inv(P^-)^{\ast}T'P^-$.
Hence $P_T$, and also $L=\Ker P_T$, continuously depend on $(\sigma,T')$.
It follows that the map $(A,T')\mapsto(A,L)$ is continuous.

Conversely, for $u\in\EmY$ we have 
$Tu = \Pmort i\sigma(n) P^+ P_{L,\,\EY^+}u$,
where $\Pmort = P^-(P^- + (P^-)^*-1)\inv$ is the orthogonal projection of $\EY$ onto $\EmY$ 
(see \eqref{eq:Port} for the formula of the orthogonal projection).
This implies
$$T' = \Pmort \; i\sigma(n) \; P^+ \; P_{L,\,\EY^+}\Pmort + (1-\Pmort).$$
Since all elements of this expression continuously depend on $(\sigma,L)$,
the map $(A,L)\mapsto(A,T')$ is continuous.
This proves the first part of the proposition.

By the definition of $T'$, we have $\chi_{(-\infty,0)}(T_x) = \chi_{(-\infty,0)}(T'_x)$,
where $\chi_S$ denotes the characteristic function of a subset $S$ of $\R$.
Hence $F_x$ considered as a point of $\Gr(E_x)$ coincides with $\im\br{\chi_{(-\infty,0)}(T'_x)}$
and thus is an analytic function of $T'_x$.
Therefore, $F$ is a smooth subbundle of $E$ and continuously depends on $T'$ in the $C^1$-topology.
Together with the continuity of $T'$ this implies 
continuity of the map $F\colon\Ellt(E)\to C^{\infty,1}(\Gr(\EY))$.
This proves the second part of the proposition.
\endproof

\section{The invariant $\Psi$ and its properties}\label{sec:F-Psi}

\upskip
\sub{Gluing of bundles.}
Let $\gamma\in\Omega_g\Ellt(E)$,
that is $\gamma\colon[0,1]\to\Ellt(E)$ is a path in $\Ellt(E)$ such that $\gamma(1)=g\gamma(0)$, $g\in U(E)$.
With every such pair $(\gamma,g)$ we associate a number of vector bundles.

First, lift $E$ to the vector bundle $\Eh = E\times [0,1]$ over $M \times [0,1]$.
Then form the vector bundle $\E$ over $\MS$ as the factor of $\Eh$,
identifying $(u,1)$ with $(gu,0)$ for every $u\in E$.

The one-parameter family $\Em_t=\Em(\gamma(t))$ of subbundles of $E$ forms the subbundle $\Eh^-$ of $\Eh$.
The condition $\gamma(1)=g\gamma(0)$ implies $\Em_1=g\Em_0$,
so $\Eh^-$ descends onto $\MS$ giving rise to the subbundle $\EEm = \EEm(\gamma,g)$ of $\E$
such that the following diagram is commutative:
\[
\begin{tikzcd}
 \Eh^- \arrow{d}\arrow[hookrightarrow]{r}
 & \Eh \arrow{d}\arrow{r}
 & M\times[0,1] \arrow{d}
\\
 \EEm \arrow[hookrightarrow]{r}
 & \E \arrow{r}
 & \MS
\end{tikzcd}
\]
In the same manner, from the one-parameter family of vector bundles $\EmY(\gamma(t))\subset\EY$ 
we construct the vector bundles $\Eh^-_{\p}\subset\Eh_{\p}$ over $\pM\times[0,1]$.
Twisting by $g$ and gluing as described above, we obtain the vector bundles $\EEmY\subset\EEY$ over $\pMS$.
Equivalently, $\EEY$ and $\EEmY$ can be obtained as the restrictions of $\E$ and $\EEm$ to $\pMS$.

The one-parameter family $F_t = F(\gamma(t))$ of subbundles of $\EmY(\gamma(t))$ 
forms the subbundle $\Fh$ of $\Eh^-_{\p}$.
Again, the condition $\gamma(1)=g\gamma(0)$ implies $F_1=gF_0$,
so $\Fh$ descends onto $\pMS$ giving rise to the subbundle $\F = \F(\gamma,g)$ of $\EEmY$
such that the following diagram is commutative:
\[
\begin{tikzcd}
   \Fh          \arrow{d}\arrow[hookrightarrow]{r}
 & \Eh^-_{\p} \arrow{d}\arrow[hookrightarrow]{r}
 & \Eh_{\p}   \arrow{d}\arrow{r}
 & \pM\times[0,1] \arrow{d}
\\
   \F    \arrow[hookrightarrow]{r}
 & \EEmY \arrow[hookrightarrow]{r}
 & \EEY  \arrow{r}
 & \pMS
\end{tikzcd}
\]
If $g=\Id$, then we will write $\F(\gamma)$ instead of $\F(\gamma,\Id)$.

\medskip

\sub{Definition of $\Psi(\gamma,g)$.}
The orientation on $M$ induces the orientation on $\pM$.
We equip $\pM$ with an orientation in such a way that the pair
(outward normal to $\pM$, positive tangent vector to $\pM$)
has a positive orientation.

The product $\pMS$ is a two-dimensional manifold, namely a disjoint union of tori.
Let $[\pMS]\in H_2(\pMS)$ be its fundamental class.
The first Chern class $c_1(\F)$ of the vector bundle $\F$ 
is an element of the second cohomology group $H^2(\pMS)$,
so one can compute its value on $[\pMS]$, 
obtaining the integer-valued invariant 
\begin{equation}\label{eq:Psi}
	\Psi(\gamma,g) = c_1(\F(\gamma,g))[\pMS].
\end{equation}
If $g=\Id$, then we will write $\Psi(\gamma)$ instead of $\Psi(\gamma,\Id)$.

\sub{The homomorphism $\psi$.}
The first Chern class is additive with respect to direct sum of vector bundles, 
so we can define the homomorphism of commutative groups $\psi\colon K^0(\pMS) \to \Z$ 
by the rule $\psi[V] = c_1(V)[\pMS]$ for any vector bundle $V$ over $\pMS$.
Then $\Psi$ can be written as 
$$\Psi(\gamma,g) = \psi[\F(\gamma,g)].$$

Consider the following three subgroups of $K^0(\pMS)$:
\begin{itemize}\upskip
	\item $G^*$ is the image of the natural homomorphism $K^0(\pM) \to K^0(\pMS)$
induced by the projection $\pMS\to\pM$.
	\item $G^{\p}$ is the image of the homomorphism $K^0(\MS) \to K^0(\pMS)$
induced by the embedding $\pMS\hookto\MS$.
	\item $G$ is the subgroup of $K^0(\pMS)$ spanned by $G^*$ and $G^{\p}$.
\end{itemize}

\begin{prop}\label{prop:psi}
The homomorphism $\psi$ is surjective with the kernel $G$. 
In other words, the following sequence is exact:
\begin{equation*}
	\begin{tikzcd}
	  0 \arrow{r} & G \arrow{r} & K^0(\pMS) \arrow{r}{\psi} & \Z \arrow{r} & 0.
	\end{tikzcd}
\end{equation*}
\end{prop}

\proof
Denote the connected components of $\pM$ by $\pM_1,\ldots,\pM_m$.
The group $K^0(\pMS)$ is isomorphic to $\Z^{2m}$, with the isomorphism given by 
$$[V]\mapsto(r_1,\ldots,r_m, a_1,\ldots,a_m),$$
where $r_j$ is the rank of the restriction $V_j$ of a vector bundle $V$ to $\pM_j\times S^1$
and $a_j = c_1(V_j)[\pM_j\times S^1]$.

In this designations, the subgroup $G^*$ consists of elements with $a_1=\ldots=a_m=0$.
The subgroup $G^{\p}$ consists of elements with $r_1=\ldots=r_m$ and $\sum_j a_j = 0$.
The span $G$ of $G^*$ and $G^{\p}$ consists of elements with $\sum_j a_j = 0$.
The homomorphism $\psi$ takes $(r_1,\ldots,r_m, a_1,\ldots,a_m)$ to $\sum_j a_j$,
so it is surjective with the kernel $G$.
This completes the proof of the proposition.
\endproof

\sub{Special subspaces.}
The following two subspaces of $\Ellt(E)$ will play special role:
\begin{itemize}\upskip
	\item $\Ellp(E)$ consists of all $(A,T)\in\Ellt(E)$ with positive definite $T$.
	\item $\Ellm(E)$ consists of all $(A,T)\in\Ellt(E)$ with negative definite $T$.
\end{itemize}

\begin{prop}\label{prop:Fpm}
Let $\gamma\in\Omega_g\Ellt(E)$. 
Then the following statements hold:
\begin{enumerate}\upskip
	\item $\F(\gamma,g)=0$ if and only if $\gamma\in\Omega_g\Ellp(E)$;
	\item $\F(\gamma,g)=\EEmY(\gamma,g)$ if and only if $\gamma\in\Omega_g\Ellm(E)$.
\end{enumerate}
\end{prop}

\proof
This follows immediately from the definition of $\F$.
\endproof

\sub{Properties of $\Psi$.}
Denote by $\Omega^*\Ellt(E)$ the subspace of $\Omega\Ellt(E)$ consisting of constant loops.

\begin{prop}\label{prop:prop-Psi}
$\Psi$ has the following properties: 
\begin{enumerate}\upskip
	\item[$(\Psi0)$] 
	            $\Psi$ vanishes on $\Omega^*\Ellt(E)$, 
               $\Omega_g\Ellp(E)$, and $\Omega_g\Ellm(E)$ for every $g\in U(E)$.
	\item[$(\Psi1)$] 
	            $\Psi$ is constant on path connected components of $\Omega_g\Ellt(E)$ for every $g\in U(E)$. 							
	\item[$(\Psi2)$] 
	            $\Psi(\gamma_0\oplus\gamma_1, g_0\oplus g_1) = \Psi(\gamma_0,g_0) + \Psi(\gamma_1,g_1)$
							 for $\gamma_i\in\Omega_{g_i}\Ellt(E_i)$, $g_i\in U(E_i)$, $i=0,1$.
\end{enumerate}
\end{prop}

\proof
$(\Psi0)$.
If $\gamma\in\Omega_g\Ellp(E)$, 
then $\F(\gamma,g)=0$, so $\Psi(\gamma,g)=0$.

If $\gamma\in\Omega_g\Ellm(E)$, 
then $\F(\gamma,g)$ is the restriction to $\pMS$ of the vector bundle $\EEm(\gamma,g)$ over $\MS$, 
so $[\F(\gamma,g)]\in G^{\p}$.

If $\gamma\in\Omega^*\Ellt(E)$, $\gamma(t)\equiv(A,L)$, 
then $\F(\gamma)$ is the lifting to $\pMS$ of the vector bundle $F(A,L)$ over $\pM$,
so $[\F(\gamma)]\in G^{*}$.

In both last cases Proposition \ref{prop:psi} implies vanishing of $\Psi$.

$(\Psi1)$.
If $\gamma_0$ and $\gamma_1$ are connected by a path $(\gamma_s)$ in $\Omega_g\Ellt(2k_M)$,
then $\F_0=\F(\gamma_0,g)$ and $\F_1=\F(\gamma_1,g)$ are homotopic via the homotopy $s\mapsto \F(\gamma_s,g)$.
It follows that the classes of $\F_0$ and $\F_1$ in $K^0(\pMS)$ coincide, and thus 
$\Psi(\gamma_0,g) = \psi[\F_0] = \psi[\F_1] = \Psi(\gamma_1,g)$.

$(\Psi2)$.
Obviously, $\F(\gamma_0\oplus\gamma_1, g_0\oplus g_1) = \F(\gamma_0,g_0)\oplus\F(\gamma_1,g_1)$.
Passing to the classes in $K^0(\pMS)$ and applying the homomorphism $\psi$, we obtain the additivity of $\Psi$.
\endproof

\section{Dirac operators}\label{sec:Dir}

\upskip

\sub{Odd Dirac operators.}
Recall that $A\in\Ell(E)$ is called a Dirac operator if
$\sigma_A(\xi)^2 = \norm{\xi}^2\Id_E$ for all $\xi\in T^{\ast}M$.
We denote by $\Dir(E)$ the subspace of $\Ell(E)$
consisting of all \textit{odd} Dirac operators,
that is, operators having the form 
\begin{equation}\label{eq-Dir-odd}
	A = \matr{ 0 & A^- \\ A^+ & 0}
	\text{ with respect to the chiral decomposition } E = \Ep(A)\oplus\Em(A).
\end{equation}
Denote by $\Dirt(E)$ the subspace of $\Ellt(E)$ consisting of all pairs $(A,L)$ such that $A\in\Dir(E)$.

The following two subspaces of $\Dirt(E)$ will play special role:
$$\Dirp(E) = \Dirt(E)\cap\Ellp(E), \quad \Dirm(E) = \Dirt(E)\cap\Ellm(E).$$

\upskip
\sub{Realization of bundles.}
In the following we will need the possibility to realize some vector bundles over $\pMS$
as $\F(\gamma)$ for some $\gamma$. 
Recall that we denoted by $k_N$ the trivial vector bundle of rank $k$ over $N$.

\begin{prop}\label{prop:Fany1}
Every smooth vector bundle $V$ over $\pM$
can be realized as $F(A,L)$ for some $k\in\N$ and $(A,L)\in\Dirt(2k_M)$.
\end{prop}

\proof
$V$ can be embed as a smooth subbundle to a trivial vector bundle $k_{\pM}$ of sufficiently large rank $k$.
Choose a smooth global field $(e_1,e_2)$ of positive oriented orthonormal frames of $TM$
and define the Dirac operator acting on sections of $k_M$ (that is, $\CC^k$-valued functions on $M$) 
by the formula $D=-i\p_1+\p_2$. 
Let $D^t$ be the operator formally adjoint to $D$.
Then 
\begin{equation}\label{eq:DDt}
  A = \matr{0 & D^t \\ D & 0 }  
\end{equation}
is an odd Dirac operator acting on sections of $k_M\oplus k_M$, and $\Em(A)=k_M$.
Let $V^{\bot}$ be the orthogonal complement of $V$ in $\EmY(A)=k_{\pM}$,
and let $L$ be the boundary condition for $A$ defined by $T=(-1)_V\oplus 1_{V^{\bot}}$.
Then $(A,L)\in\Dirt(2k_M)$ and $F(A,L)=V$,
which proves the proposition.
\endproof

\begin{prop}\label{prop:Fany2}
Every smooth vector bundle $V$ over $\pMS$
can be realized as $\F(\gamma)$ for some $k\in\N$ and $\gamma\in\Omega\Dirt(2k_M)$.
\end{prop}

\proof
$V$ can be embed as a smooth subbundle to a trivial vector bundle over $\pMS$ of sufficiently large rank $k$. 
Let $(V_t)$, $t\in S^1$, be the correspondent one-parameter family of subbundles of $k_{\pM}$.
Define the odd Dirac operator $A\in\Dir(k_M\oplus k_M)$ by formula \eqref{eq:DDt}.
Let $L_t$ be the boundary condition for $A$ 
corresponding to the automorphism $T=(-1)_{V_t}\oplus 1_{V_t^{\bot}}$ of $k_M$. 
The element $(A,L_t)\in\Dirt(2k_M)$ depends continuously on $t$, 
so the family $(A,L_t)$ defines the loop $\gamma\in\Omega\Dirt(2k_M)$.
By construction, $F(A,L_t)=V_t$, so $\F(\gamma) = V$, which completes the proof of the proposition.
\endproof

\begin{prop}\label{prop:Fany3}
Let $V$ be a smooth vector bundle over $\MS$.
Then the restriction $V_{\p}$ of $V$ to $\pMS$ 
can be realized as $\F(\gamma,g)$ for some $\gamma\in\Omega_g\Dirm(2k_M)$, $k\in\N$, $g\in U(2k_M)$.
\end{prop}

\proof
Let $k$ be the rank of $V$.
The lifting of $V$ by the map $M\times[0,1]\to\MS$ is a trivial vector bundle $k_{M\times[0,1]}$,
so we can obtain $V$ from this trivial bundle, 
gluing $k_{M\times\set{1}}$ with $k_{M\times\set{0}}$ by some unitary bundle automorphism $g\in U(k_M)$.
Let $E=k_M\oplus k_M$, $\tilde{g}=g\oplus g \in U(E)$,
and $A\in\Dir(E)$ be given by formula \eqref{eq:DDt}.
Since the symbol of $A$ is $\tilde{g}$-invariant, 
$A_1=\tilde{g}A\tilde{g}\inv$ has the same symbol as $A$,
so the path $[0,1]\ni t\mapsto A_t=(1-t)A+tA_1$ is an element of $\Omega_{\tilde{g}}\Dir(E)$.
It follows that the path $\gamma$ given by the formula $\gamma(t)=(A_t,-\Id)$ 
is an element of $\Omega_{\tilde{g}}\Dirm(E)$.
By construction, $\F(\gamma,g) = V_{\p}$.
This completes the proof of the proposition.
\endproof

\begin{prop}\label{prop:Psi-surj}
  Every integer $\lambda$ can be obtained as $\lambda=\Psi(\gamma)$
	for some $k\in\N$ and $\gamma\in\Omega\Dirt(2k_M)$.	
\end{prop}

\proof 
Every integer $\lambda$ can be obtained as the first Chern number of a smooth vector bundle over a torus.
Hence $\lambda=\psi[V]$ for some smooth vector bundle $V$ over $\pMS$.
By Proposition \ref{prop:Fany2} $V$ can be realized as $V=\F(\gamma)$ for some $\gamma\in\Omega\Dirt_M$.
We obtain $\lambda=\Psi(\gamma)$, which completes the proof of the proposition.
\endproof

\section{Universality of $\Psi$}\label{sec:uni-Psi}

\upskip
\sub{Homotopies fixing the operators.}
In this section we will deal only with such deformations of elements of $\Omega_g\Ellt(E)$ 
that fix an operator family $\A=(A_t)$ and change only boundary conditions $(L_t)$. 

Let us fix an odd Dirac operator $D\in\Dir(2_M)$.
Denote by $\delta^+\in\Omega^*\Dirp(2_M)$, resp. $\delta^-\in\Omega^*\Dirm(2_M)$ the constant loop 
taking the value $(D,\Id)$, resp. $(D,-\Id)$.
We denote by $k\delta^+$, resp. $k\delta^-$ the direct sum of $k$ copies of $\delta^+$, resp. $\delta^-$.
Notice that $\EEmY(k\delta^+) = \EEmY(k\delta^-) = \F(k\delta^-) = k_{\pMS}$ 
and $\F(k\delta^+) = 0$. 

\begin{prop}\label{prop:Fhom}
Let $\gamma\colon t\mapsto(A_t,L_t)$ and $\gamma'\colon t\mapsto(A_t,L'_t)$, $t\in[0,1]$, 
be elements of $\Omega_g\Ellt(E)$ differing only by boundary conditions.
Then the following holds.
\begin{enumerate}\upskip
	\item If $\F(\gamma,g)$ and $\F(\gamma',g)$ are homotopic subbundles of $\EEmY(\gamma,g)$,
then $\gamma$ and $\gamma'$ can be connected by a path in $\Omega_g\Ellt(E)$.
	\item If $\F(\gamma,g)$ and $\F(\gamma',g)$ are isomorphic as vector bundles,
then $\gamma\oplus k\delta^+$ and $\gamma'\oplus k\delta^+$
can be connected by a path in $\Omega_{g\oplus\Id}\Ellt(E\oplus 2k_M)$ for $k$ large enough.
	\item If $[\F(\gamma,g)] = [\F(\gamma',g)] \in K^0(\pMS)$,
then $\gamma\oplus l\delta^-\oplus k\delta^+$ and $\gamma'\oplus l\delta^-\oplus k\delta^+$
can be connected by a path in $\Omega_{g\oplus\Id\oplus\Id}\Ellt(E\oplus 2l_M \oplus 2k_M)$ 
for $l$, $k$ large enough.
\end{enumerate}
\end{prop}

\proof
Notice that $\EEmY(\gamma,g)$ depends only on operators and does not depend on boundary conditions,
so $\EEmY(\gamma,g) = \EEmY(\gamma',g)$.
Denote $\EEmY=\EEmY(\gamma,g)$, $\F = \F(\gamma,g)$, and $\F' = \F(\gamma',g)$. 

1. Let $\A=(A_t)\in\Omega_g\Ell(E)$ be the correspondent path of operators.
Denote by $\L(\A,g)$ the space of all lifts of $\A$ to $\Omega_g\Ellt(E)$.
Denote by $\L^u(\A,g)$ the subspace of $\L(\A,g)$ consisting of paths $(A_t,T_t)$ 
such that the self-adjoint automorphism $T_t$ is unitary for every $t\in[0,1]$.
The subspace $\L^u(\A,g)$ is a strong deformation retract of $\L(\A,g)$, with the retraction given by the formula
$q_s(A_t,T_t) = (A_t, (1-s+s|T_t|\inv)T_t)$.
Since $q_s$ preserves $\F$,
it is sufficient to prove the first claim of the proposition for $\gamma,\gamma'\in\L^u(\A,g)$.

For a fixed $\A$, an element $\gamma\in\L^u(\A,g)$ is uniquely defined by a subbundle $\F(\gamma,g)$ of $\EEmY(\gamma,g)$,
and every deformation of $\F$ uniquely defines the deformation of $\gamma$.
Suppose that $\F$ and $\F'$ are homotopic subbundles of $\EEmY$.
A homotopy $h_s$ between $\F$ and $\F'$ can be chosen smooth by $x\in \pM$ 
and continuous (in the $C^1$-topology) by $s,t\in [0,1]$.
As described above, such a homotopy defines a path connecting 
$\gamma$ and $\gamma'$ in $\L^u(\A,g)\subset\Omega_{g}\Ellt(E)$.
This completes the proof of the first claim of the proposition.

2. If $\F$ and $\F'$ are isomorphic as vector bundles,
then $\F\oplus 0$ and $\F'\oplus 0$ are homotopic as subbundles of $\EEmY\oplus k_{\pMS}$ for $k$ large enough.
It remains to apply the first part of the proposition to the elements
$\gamma\oplus k\delta^+$ and $\gamma'\oplus k\delta^+$ of $\Omega_{g\oplus\Id}\Ellt(E\oplus 2k_M)$.

3. The equality $[\F] = [\F']$ implies that
the vector bundles $\F$ and $\F'$ are stably isomorphic,
that is, $\F_1\oplus l_{\pMS}$ and $\F_2\oplus l_{\pMS}$ are isomorphic for some integer $l$.
It remains to apply the second part of the proposition to the elements
$\gamma\oplus l\delta^-$ and $\gamma'\oplus l\delta^-$ of $\Omega_{g\oplus\Id}\Ellt(E\oplus 2l_M)$.
\endproof

\sub{The case of different operators.}
For $\gamma\in\Omega_g\Ellt(E)$, $\gamma(t) = (A_t, T_t)$,
we denote by $\gamma^+$ the element of $\Omega_g\Ellp(E)$ given by the rule $t\mapsto (A_t,\Id)$. 

Let $\gamma_i\in\Omega_{g_i}\Ellt(E_i)$, $i=1,2$.
Consider the elements
$\gamma'_1 = \gamma_1\oplus\gamma_2^+$ and $\gamma'_2 = \gamma_1^+\oplus\gamma_2$
of $\Omega_{g_1\oplus g_2}\Ellt(E_1\oplus E_2)$. 
By Proposition \ref{prop:Fpm} $\F(\gamma'_i,g_1\oplus g_2) = \F(\gamma_i,g_i)$.
On the other hand, $\gamma'_1$ and $\gamma'_2$ differ only by boundary conditions 
and thus fall within the framework of Proposition \ref{prop:Fhom}.
In particular, from the third part of Proposition \ref{prop:Fhom} we immediately get the following.

\begin{prop}\label{prop:F-iso2}
Let $\gamma_i\in\Omega_{g_i}\Ellt(E_i)$, $i=1,2$.
Suppose that $[\F(\gamma_1,g_1)] = [\F(\gamma_2,g_2)] \in K^0(\pMS)$.
Then 
$$\gamma_1\oplus\gamma_2^+\oplus l\delta^-\oplus k\delta^+ \; \text{ and } \;
  \gamma_1^+\oplus\gamma_2\oplus l\delta^-\oplus k\delta^+$$
can be connected by a path in $\Omega_{g_1\oplus g_2\oplus\Id\oplus\Id}\Ellt(E_1\oplus E_2\oplus 2l_M \oplus 2k_M)$ 
if $l$, $k$ are large enough.
\end{prop}

\sub{Semigroup of elliptic operators.}
The disjoint union \[ \Ellt_M = \coprod_{k\in\N}\Ellt(2k_M) \]
has the natural structure of a (non-commutative) graded topological semigroup 
with respect to the direct sum of operators and boundary conditions.

The pointwise direct sum of paths defines the map 
$$\Omega_g\Ellt(2k_M)\times\Omega_{g'}\Ellt(2k'_M)\to\Omega_{g\oplus g'}\Ellt(2(k+k')_M),$$
which induces the natural structure of a (non-commutative) topological semigroup on the disjoint union 
\[ \Omega_U\Ellt_M = \coprod_{k\in\N,\, g\in U(2k_M)} \Omega_g\Ellt(2k_M) 
  = \set{(\gamma,g)\colon \gamma\in\Omega_g\Ellt(2k_M), k\in\N,\, g\in U(2k_M)}. \]
The disjoint unions 
\[ \Omega_U\Ellp_M = \coprod_{k,g}\Omega_g\Ellp(2k_M), \; 
   \Omega_U\Ellm_M = \coprod_{k,g}\Omega_g\Ellm(2k_M), \; \text{and} \;
	 \Omega^*\Ellt_M  = \coprod_{k}\Omega^*\Ellt(2k_M) \]
are subsemigroups of $\Omega_U\Ellt_M$.

\sub{Universality of $\Psi$.}
Now we are ready to state the main result of this section.

\begin{thm}\label{thm:uni-Psi}
Let $\Phi$ be a semigroup homomorphism from $\Omega_U\Ellt_M$ to a commutative monoid $\Lambda$,
which is constant on path connected components of $\Omega_U\Ellt_M$. 
Then the following two conditions are equivalent:
\begin{enumerate}\upskip
	\item $\Phi$ vanishes on $\Omega^*\Ellt_M$, $\Omega_U\Ellp_M$, and $\Omega_U\Ellm_M$.
	\item $\Phi = \vartheta\circ\Psi$ for some (unique) monoid homomorphism $\vartheta\colon \Z\to\Lambda$, that is,
	$\Phi$ has the form $\Phi(\gamma,g) = \Psi(\gamma,g)\cdot \lambda$ for some invertible constant $\lambda\in\Lambda$.
\end{enumerate}\upskip
\end{thm}
Here by \q{invertible} we mean that there is $\lambda'\in\Lambda$ inverse to $\lambda$, 
that is such that $\lambda'+\lambda=0$.

\medskip
\proof
($2\Rightarrow 1$) follows immediately from properties ($\Psi0$--$\Psi2$) of Proposition \ref{prop:prop-Psi}.

Let us prove ($1\Rightarrow 2$). 
Suppose that $\Phi$ satisfies condition (1) of the theorem.
By Proposition \ref{prop:F-iso2} the equality $[\F(\gamma_1,g_1)] = [\F(\gamma_2,g_2)]$ implies 
\begin{equation}\label{eq:uni-Psi}
	\Phi(\gamma_1\oplus\gamma_2^+\oplus l\delta^-\oplus k\delta^+, g_1\oplus g_2\oplus\Id)
		= \Phi(\gamma_1^+\oplus\gamma_2\oplus l\delta^-\oplus k\delta^+, g_1\oplus g_2\oplus\Id).
\end{equation}
Since $\Phi$ vanishes on $(\gamma_i^+,g_i)\in\Omega_U\Ellp_M$, 
$(\delta^+,\Id)\in\Omega_U\Ellp_M$, and $(\delta^-,\Id)\in\Omega_U\Ellm_M$,	
\eqref{eq:uni-Psi} implies $\Phi(\gamma_1,g_1) = \Phi(\gamma_2,g_2)$.
It follows that the homomorphism $\Phi\colon\Omega_U\Ellt_M\to\Lambda$ factors through 
the (unique) semigroup homomorphism $\varphi\colon H\to\Lambda$,
where $H$ denotes the image of $\Omega_U\Ellt_M$ in $K^0(\pMS)$:
\begin{equation*}
\begin{tikzcd}
		& \Omega_U\Ellt_M \arrow{d}{[\F]} \arrow[swap, bend right]{ddl}{\Psi} \arrow[bend left]{ddr}{\Phi} & \\
		& H \arrow[swap]{dl}{\psi} \arrow[dashed]{dr}{\varphi} & \\
		\Z \arrow[dashed]{rr}{\vartheta} & & \Lambda 
	\end{tikzcd}
\end{equation*}

Suppose that $\psi(h_1) = \psi(h_2)$ for $h_1,h_2\in H$.
By Proposition \ref{prop:psi} this implies $h_1-h_2 = \mu^*+\mu^{\p}\in K^0(\pMS)$ 
for some $\mu^*\in G^*$ and $\mu^{\p}\in G^{\p}$.

The element $\mu^{\p}$ can be written as the difference of classes 
$[j^*V_2]-[j^*V_1]$ for some (smooth) vector bundles $V_1$, $V_2$ over $\MS$, 
where $j$ denotes the embedding $\pMS\hookto\MS$.
By Proposition \ref{prop:Fany3}, $[j^*V_i]$ can be realized as $[\F(\beta_i,g'_i)]$
for some $(\beta_i,g'_i)\in\Omega_U\Dirm_M$, which gives
$\mu^{\p} = [\F(\beta_2,g'_2)] - [\F(\beta_1,g'_1)]$.

Similarly, 
by Proposition \ref{prop:Fany1} $\mu^* = [\F(\alpha_2)] - [\F(\alpha_1)]$
for some $\alpha_1,\alpha_2\in\Omega^*\Dirt_M$.

Combining all this, for liftings $(\gamma_i,g_i)$ of $h_i$ to $\Omega_U\Ellt_M$  we obtain
\[
  [\F(\gamma_1,g_1)]+[\F(\beta_1,g'_1)]+[\F(\alpha_1)] = [\F(\gamma_2,g_2)]+[\F(\beta_2,g'_2)]+[\F(\alpha_2)],
\]
Applying $\varphi$ to the both sides of this equality and taking into account that
\[
  \varphi\br{ [\F(\gamma_i,g_i)]+[\F(\beta_i,g'_i)]+[\F(\alpha_i)] }
  = \Phi(\gamma_i,g_i)+\Phi(\beta_i,g'_i)+\Phi(\alpha_i,\Id) = \Phi(\gamma_i,g_i) = \varphi(h_i),  
\]
we obtain $\varphi(h_1) = \varphi(h_2)$.
Thus the equality $\psi(h_1) = \psi(h_2)$ implies $\varphi(h_1) = \varphi(h_2)$.
On the other hand, the homomorphism $\psi\colon H\to\Z$ is surjective by Proposition \ref{prop:Psi-surj}.
It follows that $\varphi$ factors through the (unique) semigroup homomorphism $\vartheta\colon\Z\to\Lambda$.
Since $\vartheta(0) = \Phi(\Omega^*\Dirt_M) = 0$, 
$\vartheta$ is a homomorphism of monoids.

Let $\lambda=\vartheta(1)$ and $\lambda'=\vartheta(-1)$.
Then $\lambda+\lambda'=0$ and $\vartheta(n) = n\lambda$ for every $n\in\Z$.
This completes the proof of the theorem.
\endproof

\section{Deformation retraction}\label{sec:retr}

The main result of this section is Proposition \ref{prop:retr1},
where we prove that the natural embedding $\Dirt(E)\hookto\Ellt(E)$ is a homotopy equivalence.
In the rest of the paper we will need only one corollary of this result, 
namely that every element of $\Omega_g\Ellp(E)$, resp. $\Omega_g\Ellm(E)$
is connected by a path with an element of $\Omega_g\Dirp(E)$, resp. $\Omega_g\Dirm(E)$.

\sub{Sections.}
First we construct two sections, which will be used below for construction of a deformation retraction.

\begin{prop}\label{prop:p}
The map $p\colon \Ell(E)\to \Sigma(E)$ is 
surjective and has a continuous section $r\colon\Sigma(E)\to \Ell(E)$
such that $r\circ p$ is fiberwise homotopic to the identity map.
\end{prop}

\proof 
We define a section $r\colon\Sigma(E)\to\Ell(E)$ by the formula
$r(\sigma) = \br{ \sigma_1\nabla_1 + \sigma_2\nabla_2}/2 + \br{ \sigma_1\nabla_1 + \sigma_2\nabla_2}^t/2$,
where $\sigma_i = \sigma(e_i)$,
$(e_1, e_2)$ is a fixed global field of frames in $TM$,
$\nabla$ is a fixed smooth connection on $E$,
and superscript $t$ means taking of formally adjoint operator.
The operation of taking formally adjoint operator leaves invariant symbol.
Moreover, it defines a continuous transformation of the space of first order operators
with the topology defined by the inclusion to
$C^{1}\br{\End(E)}^2 \times C^{0}\br{\End(E)}$,
$\sigma_1\nabla_1 + \sigma_2\nabla_2 + a \mapsto (\sigma_1,\sigma_2,a)$.
Thus $r$ is a continuous section of $p$
and defines a trivialization of the affine bundle $\Ell(E)\to\Sigma(E)$ 
with the fiber $C^{\infty,0}\br{\End\sa(E)}$.
Thus $r\circ p$ is fiberwise homotopic to the identity map,
which completes the proof of the proposition.
\endproof

\medskip

Denote by $\Sigma^D(E) = p(\Dir(E))$ the subspace of $\Sigma(E)$ consisting of symbols of Dirac operators.

\begin{prop}\label{prop:pD}
The restriction of $p$ to $\Dir(E)$ has a continuous section 
$r^D\colon\Sigma^D(E)\to \Dir(E)$.
\end{prop}

\proof
Let $\sigma\in\Sigma^D(E)$ and $A = r(\sigma)$.
Denote by $S$ the bundle automorphism of $E$,
whose restrictions on fibers are the orthogonal reflections in the fibers of $\Em(\sigma)$.
We define $r^D(\sigma)$ by the formula $r^D(\sigma)=(A-SAS)/2$.
Obviously, it is a Dirac operator, which is odd with respect to the chiral decomposition 
$E=\Ep(\sigma)\oplus\Em(\sigma)$ and has the same symbol $\sigma$ as $A$.
Since $S$ depends continuously on $\sigma$, 
the map $r^D\colon\Sigma^D(E)\to \Dir(E)$ is a continuous section of $\restr{p}{\Dir(E)}$.
This completes the proof of the proposition.
\endproof

\sub{Retraction of symbols.}
The following proposition is the key result of this section.

\begin{prop}\label{prop:hom_Dir1}
The subspace $\Sigma^D(E)$ is a strong deformation retract of $\Sigma(E)$.
Moreover, a deformation retraction can be chosen $U(E)$-equivariant 
and preserving $\Em(\sigma)$.
\end{prop}

\proof
For any $\sigma\in\Sigma(E)$ the automorphism $Q=\sigma(e_1)\inv\sigma(e_2)$ of $E$ 
leaves the subbundles $\Em=\Em(\sigma)$ and $\Ep=\Ep(\sigma)$ invariant. 
Denote by $Q^-$ (resp. $Q^+$) the restriction of $Q$ to $\Em$ (resp. $\Ep$).
By the construction of $\Em$ and $\Ep$, 
all eigenvalues of $Q^-_x$ (resp. $Q^+_x$) have negative (resp. positive) imaginary part for every $x\in M$. 

Denote by $J$ the restriction of $\sigma(e_1)$ to $\Em$;
it is a smooth bundle isomorphism from $\Em$ onto its orthogonal complement $(\Em)^{\bot}$.

Finally, with every $\sigma\in\Sigma(E)$ we associate the quadruple
\begin{equation}\label{eq:theta}
  \vartheta(\sigma) = (\Em, \Ep, J, Q^-).	
\end{equation}
Denote by $\Theta(E)$ the set of all quadruples $(\Em, \Ep, J, Q^-)$
such that $\Em$, $\Ep$ are transversal smooth subbundles of $E$ of half rank
(that is, $\rank\Em=\rank\Ep = \frac{1}{2}\rank E$),
$J$ is a smooth bundle isomorphism of $\Em$ onto $(\Em)^{\bot}$,
and $Q^-$ is a smooth bundle automorphism of $\Em$ such that all eigenvalues of
$Q^-_x$ have negative imaginary part for every $x\in M$.

Equip $\Theta(E)$ with the topology induced by the inclusion
\begin{multline*}
\Theta(E) \hookto C^{1}\br{\Gr(E)}^2 \times C^{1}\br{\End(E)}^2, \\
(\Em, \Ep, J, Q^-) \mapsto (\Em, \Ep, J\oplus 0_{\Ep}, Q^-\oplus 0_{\Ep}).
\end{multline*}

\begin{lem}\label{lem:4}
The map \eqref{eq:theta} defines a homeomorphism between the spaces $\Sigma(E)$ and $\Theta(E)$.
\end{lem}

\proof
Let us show first that $\vartheta$ is a bijection.
Let $(\Em, \Ep, J, Q^-)\in\Theta(E)$.
Then $\sigma_1^-=J$, $\sigma_2^-=JQ^-$ are smooth bundle isomorphisms from $\Em$ onto $(\Em)^{\bot}$.

The Hermitian structure on $E$ defines the non-degenerate pairings
$\Ep_x\times(\Em_x)^{\bot}\to\CC$ and $(\Ep_x)^{\bot}\times\Em_x\to\CC$ for each $x\in M$.
Hence there exist (unique) smooth bundle isomorphisms $\sigma_1^+$, $\sigma_2^+$ from $\Ep$ onto $(\Ep)^{\bot}$
such that $\bra{\sigma_i^+ u, v} = \bra{u, \sigma_i^- v}$ for any $u\in\Ep_x$, $v\in\Em_x$, $x\in M$.
We define the endomorphism $\sigma_i$ of $E$ by the condition that
the restriction of $\sigma_i$ to $\Ep$, resp. $\Em$ coincides with $\sigma_i^+$, resp. $\sigma_i^-$.

Every elements $u,v\in E_x$ can be written as $u=u^++u^-$, $v=v^++v^-$
with $u^+, v^+\in\Ep_x$, $u^-,v^-\in\Em_x$. 
We get
$\bra{\sigma_i u, v} = \bra{\sigma_i^+ u^+, v^-} + \bra{\sigma_i^- u^-, v^+} =
   \bra{u^+, \sigma_i^- v^-} + \bra{u^-, \sigma_i^+ v^+} = \bra{u, \sigma_i v}$.
Thus $\sigma_1$ and $\sigma_2$ are self-adjoint. 

Let $(c_1,c_2)\in\R^2\setminus\set{0}$.
Then
$c_1\sigma_1^- +c_2\sigma_2^- = \sigma_1^-(c_1 +c_2Q^-)$ is an isomorphism of $\Em$ onto $(\Em)^{\bot}$. 
By definition of $\sigma_i^+$,
$\bra{(c_1\sigma_1^+ +c_2\sigma_2^+) u, v} = \bra{u, (c_1\sigma_1^- +c_2\sigma_2^-) v}$
for any $u\in\Ep_x$, $v\in\Em_x$. 
Therefore, $c_1\sigma_1^+ +c_2\sigma_2^+$ is an isomorphism of $\Ep$ onto $(\Ep)^{\bot}$.
The direct sum decompositions $\Em\oplus\Ep = E = (\Em)^{\bot}\oplus(\Ep)^{\bot}$ imply that
$c_1\sigma_1 +c_2\sigma_2$ is a smooth bundle automorphism of $E$.
Thus $(\sigma_1, \sigma_2)$ determines the self-adjoint elliptic symbol $\sigma\in\Sigma(E)$, $\sigma(e_i)=\sigma_i$.

The automorphism $Q=\sigma_1\inv\sigma_2$ of $E$ leaves $\Em$ and $\Ep$ invariant,
and the restriction of $Q$ to $\Em$ coincides with $Q^-$.
All eigenvalues of $Q^-$ have negative imaginary part.
Ranks of $\Em$ and $\Ep$ coincide,
so by Proposition \ref{prop:Q2} all eigenvalues of the restriction of $Q$ to $\Ep$ have positive imaginary part.

By construction, $\vartheta(\sigma) = (\Em, \Ep, J, Q^-)$.
The same construction shows that $\sigma$ is determined uniquely by the quadruple $(\Em, \Ep, J, Q^-)$.
Therefore $\vartheta$ defines a bijection between $\Sigma(E)$ and $\Theta(E)$.

By Proposition \ref{lem:Empcont}, $\vartheta$ is continuous.
The construction of the inverse map given above shows that $\vartheta\inv$ is also continuous.
This completes the proof of the lemma.
\endproof

\bigskip
\textbf{\textit{Continuation of the proof of Proposition \ref{prop:hom_Dir1}.}}
By this lemma, instead of a deformation retraction of $\Sigma(E)$
we can construct a deformation retraction of $\Theta(E)$ onto the subspace
$$\Theta^D(E) = \vartheta(\Sigma^D(E)) = \set{(\Em, \Ep, J, Q^-)\in\Theta(E) \colon\;
   \Ep = (\Em)^{\bot}, \; J\in U(\Em,\Ep), \; Q^-=-i\Id}.$$
For fixed $\Em$, all three ingredients of the triple $(\Ep, J, Q^-)$ can be deformed independently of one another.
We define a homotopy $h_s(\Em, \Ep, J, Q^-) = (\Em, \Ep_s, J_s, Q^-_s)$ by the formulas
$$J_s = \br{s(JJ^{\ast})^{-1/2}+1-s}J, \quad Q^-_s = -is\Id+(1-s)Q^-,$$
and $\Ep_s$ be the graph of $(1-s)B$, where
$B$ is the smooth homomorphism from $(\Em)^{\bot}$ to $\Em$ with the graph $\Ep$.

Obviously, $h_0=\Id$, the image of $h_1$ is contained in $\Theta^D(E)$,
and the restriction of $h_s$ to $\Theta^D(E)$ is the identity for all $s\in[0,1]$.
Thus $h$ defines a deformation retraction of $\Sigma(E)$ onto $\Sigma^D(E)$.
By construction, $h_s$ is $U(E)$-equivariant and preserves $\Em(\sigma)$ for every $s\in[0,1]$.
This completes the proof of the Proposition.
\endproof

\sub{Retraction of operators.}
Using results of Propositions \ref{prop:p}--\ref{prop:hom_Dir1},  
we are now able to prove the following result.

\begin{prop}\label{prop:retr1}
The natural embedding $\Dir(E)\hookto\Ell(E)$ is a homotopy equivalence.
Moreover, there exists a deformation retraction $H$
of $\Ell(E)$ onto a subspace of $\Dir(E)$
having the following properties
for all $s\in[0,1]$ and $A\in\Ell(E)$, with $A_s=H_s(A)$:
\begin{enumerate}\upskip
	\item[(1)] $\Em(A_s)=\Em(A)$.
	\item[(2)] The symbol of $A_s$ depends only on $s$ and the symbol $\sigma_A$ of $A$. 
	\item[(3)] The map $H_s\colon\sigma_{A}\mapsto\sigma_{A_s}$ defined by (2) is $U(E)$-equivariant.
	\item[(4)] If $A\in\Dir(E)$, then $\sigma_{A_s} = \sigma_{A}$. 
	\item[(5)] If $A, B\in\im H_1$ and the symbols of $A$ and $B$ coincide, then $A=B$.
\end{enumerate}
\end{prop}

We will need only properties (1-3) in this paper. 
Properties (4-5) will be used in the next paper \cite{Pr18}.

\medskip
\proof
Throughout the proof, we call a homotopy $[0,1]\times\Ell(E)\to\Ell(E)$ \q{nice}
if it satisfies conditions (1-3) of the proposition.
Obviously, the set of nice homotopies is closed under concatenation. 
We will construct a desired deformation retraction $H$ as the concatenation of three nice homotopies.
Then we show that the resulting homotopy satisfies conditions (4-5) as well. 

Let $r\colon\Sigma(E)\to \Ell(E)$ be a section from Proposition \ref{prop:p}
and $r^D\colon\Sigma^D(E)\to \Dir(E)$ be a section from Proposition \ref{prop:pD}.
The linear fiberwise homotopy $q$ between $r\circ p$ and the identity map 
is a nice deformation retraction of $\Ell(E)$ onto $r(\Sigma(E))$. 
The composition $r\circ h_s\circ p$ gives a nice deformation retraction 
of $r(\Sigma(E))$ onto $r(\Sigma^D(E))\subset p\inv(\Sigma^D(E))$;
we will denote it by the same letter $h$.
The linear fiberwise homotopy $q^D$ between $r^D\circ p$ and the identity map
is a nice deformation retraction of $p\inv(\Sigma^D(E))$ onto $r^D(\Sigma^D(E))\subset\Dir(E)$. 

\[
\begin{tikzcd}
  & & r(\Sigma^D(E)) \arrow[hookrightarrow]{d} & & \\
    \Dir(E) \arrow[hookleftarrow]{r}
	& r^D(\Sigma^D(E)) \arrow[hookrightarrow]{r} 
	& p\inv(\Sigma^D(E)) \arrow{d}[swap]{p}  \arrow[hookrightarrow]{r} \arrow[bend right]{l}[swap]{q^D_1}
	& r(\Sigma(E)) \arrow{d}[swap]{p} \arrow{lu}[swap]{h_1} \arrow[hookrightarrow]{r} 
	& \Ell(E) \arrow{ld}{p} \arrow[bend right]{l}[swap]{q_1}
	\\
  & & \Sigma^D(E) \arrow[bend right]{u}[swap]{r} \arrow{lu}{r^D} \arrow[hookrightarrow]{r} 
	& \Sigma(E) \arrow[bend right]{u}[swap]{r} \arrow[bend left]{l}{h_1} &
\end{tikzcd}
\]
Concatenating $q$, $h$, and $q^D$,
we obtain a nice deformation retraction $H$ of $\Ell(E)$ onto the subspace $r^D(\Sigma^D(E))$ of $\Dir(E)$:
$$H_s(A) = \case{q_{3s}(A) \text{ for } 0\leq s\leq 1/3, \\ 
                   h_{3s-1}q_1(A) \text{ for } 1/3\leq s \leq 2/3 \\
									 q^D_{3s-2}h_1 q_1(A) \text{ for } 2/3\leq s \leq 1. }$$

If $A\in\Dir(E)$, then $\sigma_A\in\Sigma^D(E)$, so the symbol of $A_s$ is independent of $s$.
If $A\in\im H_1$, then $A = r^D(\sigma_A)$.
This proves conditions (4-5) of the proposition.

It remains to check that the natural embedding $\Dir(E)\hookto\Ell(E)$ is a homotopy equivalence.
For every $A\in\Dir(E)$, the image $H_1(A) = A_1$ also lies in $\Dir(E)$, 
but we need to be careful because $A_s$ is not necessarily odd for $s\in(0,1)$.
By property (4) the symbols of $A_1$ and $A$ coincide. 
Thus the formula $H'_s(A) = (1-s)A+sH_1(A)$ defines a continuous map
$H'\colon[0,1]\times\Dir(E)\to\Dir(E)$ such that $H'_0=\Id$ and $H'_1 = H_1$.
It follows that the restriction of $H_1$ to $\Dir(E)$ and the identity map $\Id_{\Dir(E)}$
are homotopic as maps from $\Dir(E)$ to $\Dir(E)$.
On the other hand, the map $H_1\colon\Ell(E)\to\Ell(E)$ is homotopic to $\Id_{\Ell(E)}$ via the homotopy $H$.
It follows that $H_1\colon\Ell(E)\to\Dir(E)$ is homotopy inverse to the embedding $\Dir(E)\hookto\Ell(E)$,
that is this embedding is a homotopy equivalence.									
This completes the proof of the proposition.
\endproof

\begin{prop}\label{prop:retr1t}
The natural embedding $\Dirt(E)\hookto\Ellt(E)$ is a homotopy equivalence.
Moreover, there exists a deformation retraction of $\Ellt(E)$ onto a subspace of $\Dirt(E)$
preserving both $\Em(A)$ and $F(A,L)$.
\end{prop}

\proof
Since the deformation retraction $H$ constructed in Proposition \ref{prop:retr1} preserves $\Em(A)$,
one can define the deformation retraction $\bar{H}\colon[0,1]\times\Ellt(E)\to\Ellt(E)$ covering $H$
and satisfying the conditions of the proposition 
by the formula $\bar{H}_s(A,T)=(H_s(A),T)$ for $(A,T)\in\Ellt(E)$.
\endproof

\medskip
\sub{Retraction of paths.}
Applying the deformation retraction from last two propositions pointwise 
and slightly correcting it on the ends of a path,
we obtain a deformation retraction of the space of paths in $\Ell(E)$ and in $\Ellt(E)$.

\begin{prop}\label{prop:retr2}
Let $g\in U(E)$. Then the following holds.
\begin{enumerate}\upskip
	\item There exists a deformation retraction of $\Omega_g\Ell(E)$ onto a subspace of $\Omega_g\Dir(E)$
preserving $\EEm(\gamma,g)$ for every $\gamma\in\Omega_g\Ell(E)$.
	\item There exists a deformation retraction of $\Omega_g\Ellt(E)$ onto a subspace of $\Omega_g\Dirt(E)$
preserving both $\EEm(\gamma,g)$ and $\F(\gamma,g)$.
\end{enumerate}
\end{prop}

\proof
1. Let $\rho_0,\rho_1\colon[0,1]\to\R$ be a	partition of unity subordinated to the covering 
$[0,1] = U_0\cup U_1$, $U_0=[0,2/3)$, $U_1=(1/3,1]$,
that is, $\supp\rho_i\subset U_i$ and $\rho_0+\rho_1\equiv 1$.
Let $h\colon[0,1]\times\Ell(E)\to\Ell(E)$ be a deformation retraction 
of $\Ell(E)$ onto a subspace of $\Dir(E)$ satisfying conditions of Proposition \ref{prop:retr1}. 
Then a desired deformation retraction $[0,1]\times\Omega_g\Ell(E)\to\Omega_g\Ell(E)$ 
can be defined by the formula 
\begin{equation}\label{eq:As}
	(s,\A)\mapsto \A_s = \rho_0\A^0_{s}+\rho_1\A^1_{s},
  \text{ where } A^0_{s}(t) = h_s(\A(t)), \; \A^1_{s}(t) = gh_s(g\inv\A(t)).
\end{equation}
Indeed, by property (3) of Proposition \ref{prop:retr1} 
the operators $\A^0_{s}(t)$ and $\A^1_{s}(t)$ have the same symbols, 
so their convex combination $A_s(t)$ lies in $\Ell(E)$ for every $t\in[0,1]$.
The symbols and the chiral decompositions of the odd Dirac operators $\A^0_1(t)$ and $\A^1_1(t)$ coincide,
so their convex combination $\A_1(t)$ lies in $\Dir(E)$.
For $s=0$ we get $\A^0_0 = \A^1_0 = \A$, so $\A_0 = (\rho_0+\rho_1)\A = \A$.
For each $s\in[0,1]$ we have 
$$\A_s(1) = \A^1_{s}(1) = gh_s(g\inv\A(1)) = gh_s(\A(0)) = g\A^0_{s}(0) = g\A_s(0),$$
so $\A_s$ lies in $\Omega_g\Ell(E)$.
Since $\Em(A)$ depends only on the symbol of $A$ and is preserved by $h_s$, 
we get $\EEm(\A_s,g) = \EEm(\A,g)$ for every $s\in[0,1]$.

2. We define the deformation retraction $H\colon[0,1]\times\Omega_g\Ellt(E)\to\Omega_g\Ellt(E)$
by the formula $H_s(\gamma)(t)=(\A_s(t),T(t))$ for $\gamma\in\Omega_g\Ellt(E)$,
where $\A$ is the projection of $\gamma$ to $\Omega_g\Ell(E)$,
$\gamma(t)=(\A(t),T(t))$, and $\A_s$ is defined by the formula \eqref{eq:As}.
Since $\EEm(\A_s,g)=\EEm(\A,g)$, $H_s(\gamma)$ is correctly defined.
\endproof

\sub{Deformation retraction of special subspaces.}
Let $\Ellp(E)$, resp. $\Ellm(E)$ be the subspace of $\Ellt(E)$
consisting of all $(A,L)$ with positive definite $T$, resp. negative definite $T$
(see Proposition \ref{prop:T}).

\begin{prop}\label{prop:def-pm}
For every $g\in U(E)$,
there exists a deformation retraction of $\Omega_g\Ellp(E)$ onto a subspace of $\Omega_g\Dirp(E)$
and a deformation retraction of $\Omega_g\Ellm(E)$ onto a subspace of $\Omega_g\Dirm(E)$.
\end{prop}

\proof
Let $H$ be a deformation retraction of $\Omega_g\Ellt(E)$ onto a subspace of $\Omega_g\Dirt(E)$
satisfying conditions of Proposition \ref{prop:retr2}.

For $\gamma\in\Omega_g\Ellp(E)$ and $\gamma_s = H_s(\gamma)$ we have $\F(\gamma_s)=\F(\gamma) = 0$, 
so by Proposition \ref{prop:Fpm} $\gamma_s\in\Omega_g\Ellp(E)$ for every $s$.
In particular, $\gamma_1\in\Omega_g\Ellp(E)\cap\Omega_g\Dirt(E) = \Omega_g\Dirp(E)$.

For $\gamma\in\Omega_g\Ellm(E)$ and $\gamma_s = H_s(\gamma)$ we have 
$\F(\gamma_s) = \F(\gamma) = \EEm(\gamma)=\EEm(\gamma_s)$,
so by Proposition \ref{prop:Fpm} $\gamma_s\in\Omega_g\Ellm(E)$ for every $s$.
In particular, $\gamma_1\in\Omega_g\Ellm(E)\cap\Omega_g\Dirt(E) = \Omega_g\Dirm(E)$.
\endproof

\section{Vanishing of the spectral flow}\label{sec:sf-0}

\upskip

\sub{Invertible Dirac operators.}
We have no means to detect the invertibility of an arbitrary element of $\Ellt(E)$ by purely topological methods.
However, there is a big class of \textit{odd Dirac} operators which are necessarily invertible.

\begin{prop}\label{prop:zero1}
Let $A\in\Dir(E)$, that is, $A$ is an odd Dirac operator. 
Let $T$ be a positive definite automorphism of $\EmY(A)$, 
and let $L$ be the boundary condition for $A$ defined by \eqref{eq:ortL}.
Then $A_L$ has no zero eigenvalues.
The same is true for negative definite $T$. 
In other words, both $\Dirp(E)$ and $\Dirm(E)$ are subspaces of $\Ellt^0(E)$.
\end{prop}

\upskip
This proposition explains why we distinguish odd Dirac operators.
If $T$ is definite, but $A$ is not odd, then $A_L$ no longer has to be invertible.

\medskip
\proof
Let $A$ be defined by formula \eqref{eq-Dir-odd}.
Denote the symbol of $A^+$ by $\sigma^+$.
Let $u = (u^+,u^-)$ be a section of the vector bundle $E = \Ep(A)\oplus\Em(A)$.
If $u\in\dom(A_L)$, then the restriction of $u$ to $\pM$ satisfies $i\sigma^+(n)u^+ = Tu^-$.
Since $A^+$ and $A^-$ are formally conjugate one to another,
Green's formula gives
\begin{multline*}
\int_{\pM} \bra{ Tu^-, u^- }\d l =
\int_{\pM} \bra{ i\sigma^+(n)u^+, u^- }\d l
 = \int_M \br{ \bra{ A^+ u^+, u^- } - \bra{ u^+, A^- u^- } } \d s,
\end{multline*}
where $\d l$ is the length element on $\pM$ and $\d s$  is the volume element on $M$.

Suppose now that $A_L u=0$.
Then $A^+ u^+ = A^- u^- = 0$, so the last integral vanishes and we obtain
$\int_{\pM} \bra{ Tu^-, u^- }\d l = 0$.
If $T$ is positive definite or negative definite on $\pM$,
then the last equality implies vanishing of $u^-$ on $\pM$.
This together with the boundary condition $i\sigma^+(n)u^+ = Tu^-$ 
implies vanishing of $u^+$  on $\pM$. 
By the weak inner unique continuation property of Dirac operators \cite{BBW-93},
we get $u\equiv 0$ on whole $M$.
It follows that $A_L$ has no zero eigenvalues,
which completes the proof of the proposition.
\endproof

\sub{Vanishing of the spectral flow.}
Our next goal is to show that the spectral flow satisfies the first condition of Theorem \ref{thm:uni-Psi}.

\begin{prop}\label{prop:sf0}
Let $\gamma$ be an element of $\Omega^*\Ellt_M$, $\Omega_U\Ellp_M$, or $\Omega_U\Ellm_M$.
Then $\gamma$ is connected by a path with an element of $\Omega_U\Ellt^0_M$, 
and hence $\spf(\gamma) = 0$.
\end{prop}

\proof
Suppose that $\gamma\in\Omega_g\Ellp(E)$ or $\Omega_g\Ellm(E)$, $g\in U(E)$.
By Proposition \ref{prop:def-pm}, $\gamma$ is connected by a path with an element $\gamma_1$
of $\Omega_g\Dirp(E)$ or $\Omega_g\Dirm(E)$ respectively.
By Proposition \ref{prop:zero1} $\gamma_1\in\Omega_g\Ellt^0(E)$.

Suppose that $\gamma\in\Omega^*\Ell(E)$, that is, $\gamma(t)\equiv(A,L)$.
Since $A_L$ is Fredholm, $A_L-\lambda$ is invertible for some $\lambda\in\R$.
The path $\gamma_s(t) = (A-s\lambda,L)$ connects $\gamma$ with
the constant loop $\gamma_1\in\Omega^*\Ell^0(E)$.
	
Since the spectral flow vanishes on paths in $\Ellt^0(E)$, $\spf(\gamma_1) = 0$.
The homotopy invariance of the spectral flow implies $\spf(\gamma) = 0$,
which completes the proof of the proposition.
\endproof

\section{The spectral flow formula}\label{sec:sf-proof}\label{sec:proof1}

Now we are ready to compute the spectral flow.

\begin{thm}\label{thm:sf}
Let $\gamma\colon[0,1]\to\Ellt(E)$ be a continuous path such that $\gamma(1)=g\gamma(0)$
for some smooth unitary bundle automorphism $g$ of $E$.
Then
$\sp{\gamma} = \Psi(\gamma,g)$.
\end{thm}

The proof consists of a sequence of lemmas.

\begin{lem}\label{lem:sf-lambda}
   There is an integer $\lambda=\lambda_M$ depending only on $M$ such that
\begin{equation}\label{eq:sf-const}
  \sp{\gamma} = \Psi(\gamma,g)\cdot\lambda
\end{equation}
for every $\gamma\in\Omega_g\Ellt(E)$, $g\in U(E)$.
\end{lem}

\proof
The spectral flow defines the homomorphism $\spf\colon\Omega_U\Ellt_M\to\Z$, $(\gamma,g)\mapsto\spf(\gamma)$,
which is constant on path connected components of $\Omega_U\Ellt_M$.
By Proposition \ref{prop:sf0}, the spectral flow vanishes on 
$\Omega^*\Ellt_M$, $\Omega_U\Ellp_M$, and $\Omega_U\Ellm_M$.
Thus $\Phi=\spf$ and $\Lambda=\Z$ satisfy the first condition of Theorem \ref{thm:uni-Psi}. 
By Theorem \ref{thm:uni-Psi} there is a $\lambda\in\Z$ such that \eqref{eq:sf-const} holds
for every $\gamma\in\Omega_g\Ellt(2k_M)$.
Since every vector bundle over $M$ is trivial, this completes the proof of the lemma.
\endproof

\begin{lem}\label{lem:metric}
   The value of $\lambda$ does not depend on the choice of a metric on $M$.
\end{lem}

\proof
Let $h$, $h'$ be two metrics on $M$.
The Hilbert spaces $L^2(M,h;E)$ and $L^2(M,h';E)$ are isomorphic,
with an isometry given by the formula $u \mapsto cu$,
where $c$ is the positive-valued function on $M$ defined by the formula $c = \sqrt{\det(h')/\det(h)}$. 
This isometry induces the bijection between the spaces $\Ellt(M,h;E)$ and $\Ellt(M,h';E)$
and leaves invariant the spectral flow of paths.
On the other hand, such an isometry leaves invariant both the symbols of operators and local boundary conditions,
so it leaves invariant $F(A,L)$.
The conjugation by $c$ also leaves invariant bundle automorphism $g$.
Therefore, the aforementioned bijection $\Ellt(M,h;E)\to\Ellt(M,{h'};E)$
does not affect $\F(\gamma,g)$.
This implies that the factor $\lambda$ in \eqref{eq:sf-const} is the same for metrics $h$ and $h'$.
Since $h$ and $h'$ are arbitrary metrics, $\lambda$ does not depend on the choice of a metric.
\endproof

\begin{lem}\label{lem:ann}
  If $M$ is diffeomorphic to the annulus, then $\lambda_M = \lambda\ann =1$.
\end{lem}

\proof
This was proven by the author in \cite[Theorem 4]{Pr} ($\lambda\ann$ is denoted by $c_2$ there).
The proof is based on the direct computation of the spectral flow for
the Dirac operator on $S^1\times[0,1]$ with varying connection and 
fixed boundary condition.
\endproof

\begin{lem}\label{lem:lambda}
  For any smooth oriented connected surface $M$ the values of $\lambda_M$ and $\lambda\ann$ coincide.
\end{lem}

\proof 
There are different ways to reduce the computation of $\lambda_M$ 
to the case of an annulus.
Here we describe one of them,
namely the splitting of $M$ into two pieces: the smaller surface $M'$ diffeomorphic to $M$
and the collar $M''$ of the boundary.
Following ideas of P.~Kirk and M.~Lesch from \cite{KL}, 
we take the Dirac operator which has the product form near boundary
and choose mutually orthogonal boundary conditions on the sides of the cut.
Then the spectral flow over $M$ coincides with the sum of spectral flows over $M'$ and $M''$.
Since $M''$ is the disjoint union of annuli,
this reasoning allows to reduce the computation of $\lambda_M$ to the computation for the annuli.
Let us describe this procedure in more detail.

Let $U$ be a collar neighbourhood of $\pM$ in $M$;
we identify $U$ with the product $(-2\varepsilon,0]\times\pM$. 
Let $(y,z)$ be the coordinates on $U$, with $y\in \pM$, $z\in(-2\varepsilon,0]$,
and $(\p_z,\p_y)$ a positive oriented basis in $TU$.
Equip $M$ with a metric whose restriction to $U$ has the product form $dl^2 = dy^2+dz^2$.

Let $D\in\Dir(E)$ be an odd Dirac operator acting on sections of 
$E=E^+(D)\oplus E^-(D)$ with $E^+(D) = E^-(D) = 2_M$.
Adding a bundle automorphism to $D$ if required, 
we can ensure that the restriction of $D$ to $U$ 
has the product form $\restr{D}{U}=-i(\sigma_1\p_z + \sigma_2\p_y)$,
where $\sigma_1 = \smatr{0 & \sigma_1^-  \\ \sigma_1^+ & 0}$, 
$\sigma_2 = \sigma_1 Q$, $Q=\smatr{i & 0 \\ 0 & -i}$.

Let $\F$ be a vector bundle of rank $1$ over $\pMS$ such that $c_1(\F)[\pMS] \neq 0$.
Choose the smooth embedding of $\F$ into the trivial vector bundle of rank $2$ over $\pMS$.
Restricting this embedding to $\pM\times\set{t}$,
we obtain the smooth loop $(F_t)_{t\in S^1}$ of smooth subbundles $F_t$ of $2_{\pM}$.
Define the smooth automorphisms $T_t$ of $2_{\pM}$
by the formula $T_t = (-1)_{F_t}\oplus 1_{F^{\bot}_t}$.
Let $L_t\subset \EY$ be the correspondent boundary condition for $D$
(that is, $L_t$ is obtained from $T_t$ as described in Proposition \ref{prop:T}).
Then $\F=\F(\gamma)$ for the loop $\gamma\in\Omega\Ellt(E)$ defined by the formula $\gamma(t)=(D,L_t)$.
By Lemma \ref{lem:sf-lambda},
$$\sp{D,L_t} = c_1(\F)[\pMS]\cdot \lambda_M.$$

Let us cut $M$ along $N = \set{-\varepsilon}\times\pM \subset U$.
We obtain the disconnected surface $M\cut = M'\amalg M''$,
where $M'' = [-\varepsilon,0]\times\pM$ is the disjoint union of annuli
and $M'$ is diffeomorphic to $M$.
Denote by $E\cut = E'\amalg E''$ the lifting of $E$ on $M\cut$,
and  by $D\cut = D'\amalg D''$ the lifting of $D$ on $M\cut$.
By $N'$, $N''$ denote the sides of the cut, so that $\pM'=N'$ and $\pM''=N''\amalg \pM$.

The restriction of $E\cut$ to $N'\amalg N''$ is isomorphic 
to the disjoint union of two copies of $\restr{E}{N}$.
Let us identify its sections with sections of the vector bundle 
$\bar{E}_{\p} = \restr{(E\oplus E)}{N}$. 
The diagonal subbundle $\Delta = \set{u\oplus u}$ of $\bar{E}_{\p}$ 
defines the so called transmission boundary condition on the cut.
The natural isometry $L^2(E) \to L^2(E\cut)$ takes the operator
$D_{L_t}$ to the operator $D\cut_{\Delta\amalg L_t}$.
Therefore, $D\cut_{\Delta\amalg L_t}$ is a self-adjoint Fredholm regular operator on $L^2(E\cut)$, and
$$\sp{D,L_t} = \sp{D\cut, \Delta\amalg L_t}.$$
Extending the identification above to the identification of sections of
$\restr{E\cut}{U'\amalg U''}$ with sections of $\bar{E} = \restr{(E\oplus E)}{U'}$, 
where $U'=(-2\varepsilon,-\varepsilon]\times\pM$, $U''=[-\varepsilon,0)\times\pM$,  
we can write $D\cut$ in the collar of the cut as
\[
 \bar{D} = -i(\bar{\sigma}_1\p_{\bar{z}}+\bar{\sigma}_2\p_y),
  \text{ where } \bar{\sigma}_1=\matr{\sigma_1 & 0 \\ 0 & -\sigma_1}, \;
	\bar{\sigma}_2=\matr{\sigma_2 & 0 \\ 0 & \sigma_2},
\]
and $\bar{z}$ is the normal coordinate increasing in the direction of the cut 
(so $\bar{z} = z$ on $U'$ and $\bar{z} = -z-2\varepsilon$ on $U''$). 
We also change the orientation on $M'$, so that $(\p_{\bar{z}}, \p_y)$ becomes a negative oriented basis.
Then $\bar{Q} = -\bar{\sigma}_1\inv\bar{\sigma}_2 = (-Q)\oplus Q$ and
\begin{equation}\label{eq:barE}
  \bar{E}^+ = E'^-\oplus E''^+, \quad \bar{E}^- = E'^+\oplus E''^-.	
\end{equation}
The restriction $\bar{\sigma}_1^+$  of $\bar{\sigma}_1$ to $\bar{E}^+$ has the form
$\bar{\sigma}_1^+ = \sigma_1^-\oplus(-\sigma_1^+)$ with respect to decompositions \eqref{eq:barE}.

Proposition \ref{prop:T} associates with every self-adjoint automorphism $\bar{T}$ of $\bar{E}^-_{\p}$
the subbundle $\bar{L}(\bar{T})$ of $\bar{E}_{\p}$ given by the formula 
$i\bar{\sigma}_1^+\bar{u}^+ = \bar{T}\bar{u}^-$.
Each $\bar{L}(\bar{T})$ is a self-adjoint well posed boundary condition for $\bar{D}$ on the cut,
so $\bar{L}(\bar{T})\amalg L_t$ is a self-adjoint well posed boundary condition for $D\cut$.

The transmission boundary condition $\Delta$ corresponds to the unitary self-adjoint automorphism
$$\bar{T}_{\Delta} = i\bar{\sigma}_1^+ \matr{0 & 1  \\ 1 & 0} = \matr{0 & i\sigma_1^-  \\ -i\sigma_1^+ & 0}$$ 
of $\bar{E}^-_{\p} = 4_{\pM}$.
Over every point $x\in\pM$ the trace of $\bar{T}_{\Delta}$ is zero,
so it has exactly two positive and two negative eigenvalues.
$\bar{T}_{\Delta}$ can be identified with a map from $\pM$ to the complex Grassmanian $\Gr(2,4)$.
Since $\Gr(2,4)$ is simply connected,
every two maps from $\pM$ to $\Gr(2,4)$ are homotopic. 
Thus $\bar{T}_0 = \bar{T}_{\Delta}$ can be connected by a smooth homotopy 
$(\bar{T}_s)$ with $\bar{T}_1 = (-1)\oplus 1$ 
in the space of (unitary) self-adjoint bundle automorphisms of $\bar{E}^-_{\p}$.

Denote by $\bar{L}_s$ the subbundle of $\bar{E}_{\p}$ corresponding to $\bar{T}_s$, and let $\bar{L} = \bar{L}_1$.
Then $\bar{L}_s \amalg L_t$ is a self-adjoint well posed global boundary condition for $D\cut$,
so $D\cut_{\bar{L}_s\amalg L_t}$ is a regular self-adjoint Fredholm operator on $L^2(E\cut)$ for each $s,t$.
By Lemma \ref{lem:local} from the Appendix, the map
$$[0,1]\times S^1 \to \Gr\br{ H^{1/2}(\bar{E}_{\p}) \oplus H^{1/2}(\EY)} \cong \Gr\br{H^{1/2}(\EY\cut)},$$
$(s,t) \mapsto H^{1/2}(\bar{L}_s)\oplus H^{1/2}(L_t)$, is continuous.
By Proposition \ref{prop:trace}, this implies the continuity of the map
$$[0,1]\times S^1 \to \FR\sa(L^2(E\cut)) , \quad (s,t) \mapsto D\cut_{\bar{L}_s\amalg L_t}.$$
Therefore, by the homotopy invariance property of the spectral flow we have
\[
  \sp{D\cut, \Delta\amalg L_t} = \sp{D\cut, \bar{L}\amalg L_t}.
\]
The boundary condition $\bar{L}$ is given by the formula $i\bar{\sigma}_1^+\bar{u}^+ = \bar{T}_1\bar{u}^-$. 
Coming back from $\bar{E}_{\p}$ to $\restr{E\cut}{N'\amalg N''}$,
we obtain $L'\amalg L''$ in place of $\bar{L}$,
where $L'$ is the subbundle of $\restr{E\cut}{N'}$ given by the formula $i(-\sigma_{1}^+)u'^+=u'^-$ 
and $L''$ is the subbundle of $\restr{E\cut}{N''}$ given by the formula $i\sigma_{1}^+u''^+=u''^-$.
Therefore, $\bar{L}\amalg L_t$ is a \textit{local} boundary condition for $D\cut$.
Applying Lemma \ref{lem:sf-lambda} to the connected components of $M\cut$,
we obtain
\begin{multline*}
\sp{D\cut,\bar{L}\amalg L_t} = \sp{D',L'}+\sp{D'', L''\amalg L_t} = \sp{D'', L''\amalg L_t} = \\
	    \br{ c_1(\F'')[N''\times S^1] + c_1(\F)[\pMS] } \cdot \lambda\ann = c_1(\F)[\pMS]\cdot \lambda\ann,
\end{multline*}
since $\F''$ is zero vector bundle.

Combining all this together, we obtain
\[
  c_1(\F)[\pMS]\cdot \lambda_M = \sp{D,L_t} =
	\sp{D\cut, L'\amalg L''\amalg L_t}
    = c_1(\F)[\pMS]\cdot \lambda\ann.
\]
The value of $c_1(\F)[\pMS]$ does not vanish due to the choice of $\F$.
Therefore, $\lambda\ann = \lambda_M$,
which completes the proof of the lemma and of Theorem \ref{thm:sf}.
\endproof

\section{Universality of the spectral flow}\label{sec:uni-sf}

The direct sum of two invertible operators is again invertible, so
the disjoint union \[ \Ellt^0_M = \coprod_{k\in\N}\Ellt^0(2k_M) \]
is a subsemigroup of $\Ellt_M$.
The disjoint union 
\[ \Omega_U\Ellt^0_M = \coprod_{k\in\N,\, g\in U(2k_M)} \Omega_g\Ellt^0(2k_M). \]
is a subsemigroup of $\Omega_U\Ellt_M$.

\begin{thm}\label{thm:uni-sf}
Let $\Phi$ be a semigroup homomorphism from $\Omega_U\Ellt_M$ to a commutative monoid $\Lambda$,
which is constant on path connected components of $\Omega_U\Ellt_M$. 
Then the following two conditions are equivalent:
\begin{enumerate}\upskip
	\item $\Phi$ vanishes on $\Omega_U\Ellt^0_M$.
	\item $\Phi = \vartheta\circ\spf$ for some (unique) monoid homomorphism $\vartheta\colon \Z\to\Lambda$,
	      that is, $\Phi$ has the form $\Phi(\gamma,g) = \spf(\gamma)\cdot \lambda$ 
				for some invertible constant $\lambda\in\Lambda$.  
\end{enumerate}\upskip
\end{thm}

\upskip
In other words, the spectral flow defines an isomorphism of monoids
$$\spf\colon\pi_0(\Omega_U\Ellt_M)/\pi_0(\Omega_U\Ellt^0_M) \to\Z.$$ 

\proof
($2\Rightarrow 1$) follows immediately from properties (S0-S2) of the spectral flow, see Section \ref{sec:sfR}.

Let us prove ($1\Rightarrow 2$). 
Suppose that $\Phi$ satisfies condition (1) of the theorem.
By Proposition \ref{prop:sf0} $\Phi$ vanishes on $\Omega^*\Ellt_M$, $\Omega_U\Ellp_M$, and $\Omega_U\Ellm_M$.
Theorem \ref{thm:uni-Psi} then implies $\Phi = \vartheta\circ\Psi$ 
for some (unique) monoid homomorphism $\vartheta\colon \Z\to\Lambda$.  
By Theorem \ref{thm:sf} $\Psi$ is equal to the spectral flow.
Taking this all together, we obtain $\Phi = \vartheta\circ\spf$.
Taking $\lambda=\vartheta(1)$, we obtain $\Phi(\gamma,g) = \spf(\gamma)\cdot \lambda$,
which completes the proof of the theorem.
\endproof

\appendix 

\section{Appendix. Criteria of graph continuity}

\renewcommand{\thesection}{\Alph{section}} 

In this Appendix we give some general conditions describing if a family of closed operators
(in particular, differential operators on a manifold with boundary) is graph continuous.
We use the results of subsection \ref{sec:app-diff}
in the main part of the paper for two purposes.
First, Proposition \ref{prop:AL2} arises as a particular case of Proposition \ref{prop:AL2-app}.
Second, Proposition \ref{prop:trace} and Lemma \ref{lem:local} 
provide the continuity of the family of global boundary value problems
used in the proof of Lemma \ref{lem:lambda}.

After the main part of the Appendix (namely, the case of Hilbert spaces) was written, 
the author discovered that some of the results of subsections \ref{app:inj} and \ref{app:cl}, 
though in a different form and with different proofs, 
are contained in the Appendix to the recent paper of Booss-Bavnbek and Zhu \cite{BBZ-14}.
In particular, our Proposition \ref{prop:Grj} is a corollary of \cite[Proposition A.6.2]{BBZ-14}
and our Proposition \ref{prop:bounded} is a special case of \cite[Corollary A.6.4]{BBZ-14}.
Nevertheless, we leave these results and their proofs in the paper 
for the sake of completeness, and also because their statements better meet our needs.
For Hilbert spaces, our proofs have the advantage of not using elaborated estimates and inequalities.
We also add the more general case of Banach spaces to the Appendix
with the purpose of better matching the results of \cite{BBZ-14},
though we use only Hilbert spaces in the main part of the paper.

It is worth noticing that our Proposition \ref{prop:q-top} gives 
an equivalent definition of the gap topology on the space $\Gr(H)$ 
of all complemented closed linear subspaces of a Banach space $H$.
Namely, the gap topology on $\Gr(H)$ coincides with the quotient topology 
induced by the map $\Proj(H)\to\Gr(H)$, $P\mapsto\im P$,
where $\Proj(H)$ is the space of all idempotents in $\B(H)$ with the norm topology.
The author does not know if this fact was noted before.

\subsection{Complementary pairs of subspaces} 

\sub{Subspaces of a Banach space.}
Let $H$ be a Banach space. 
Denote by $\BB(H)$ the space of all bounded linear operators on $H$ with the norm topology.
Denote by $\Proj(H)$ the subspace of $\B(H)$ consisting of all idempotents.

A closed subspace $L\subset H$ is called complemented if 
there is another closed subspace $M\subset H$ such that $L\cap M=0$, $L+M=H$;
such pair $(L,M)$ is called a complementary pair.
Equivalently, $L\subset H$ is complemented if it is the image of some $P\in\Proj(H)$;
$(L,M)$ is a complementary pair if it is equal to $(\im P, \Ker P)$ for some $P\in\Proj(H)$.

We denote by $\Gr(H)$ the set of all complemented closed linear subspaces of $H$,
and by $\Grr(H)$ the set of all complementary pairs of subspaces of $H$.
We will also write $\Gr^2(H)$ instead of $\Gr(H)^2$ for convenience.

For $(L,M)\in\Grr(H)$ we denote by $P_{L,M}$ the projection of $H$ onto $L$ along $M$.
For $M\in\Gr(H)$ denote by $\Gr^M(H) = \set{ L\in\Gr(H) \colon (L,M)\in\Grr(H) }$
the set of all complement subspaces for $M$.

\begin{prop}\label{prop:PQST}
  Let $H$ be a Banach space and $P,Q\in\Proj(H)$. Then the following two conditions are equivalent:
	\begin{enumerate}\upskip
		\item Both $(\im P, \im Q)$ and $(\Ker P, \Ker Q)$ are complementary pairs of subspaces.
		\item $P-Q$ is invertible.
	\end{enumerate}\upskip
	If this is the case, then for the projection $S$ on $\im P$ along $\im Q$ 
	and the projection $T$ on $\Ker P$ along $\Ker Q$ we have:	
	\begin{equation}\label{eq:PQST}
		S=P(P-Q)\inv, \quad T=(P-1)(P-Q)\inv, \quad (P-Q)\inv = S-T, 
	\end{equation}
	and $P+Q=(2S-1)(P-Q)$ is also invertible.
\end{prop}

\proof
$(1\Rightarrow 2)$ 
Let $(\im P, \im Q), (\Ker P, \Ker Q) \in\Grr(H)$.
Denote by $S$, $T$ the elements of $\Proj(H)$ corresponding these two pairs of complementary subspaces.
Using the identities $SP=P$, $TQ=T$, $SQ=0$, and $(1-T)(1-P)=0$, we obtain
$$(S-T)(P-Q) = T+P-TP = 1-(1-T)(1-P) = 1.$$
Similarly, we have
$$(P-Q)(S-T) = Q+S-QS = 1-(1-Q)(1-S) = 1.$$
Therefore, $P-Q$ is invertible with $S-T$ the inverse operator.

$(2\Rightarrow 1)$ 
Let $P-Q$ be invertible.
It vanishes on the intersections $\im P\cap\im Q$ and $\Ker P\cap\Ker Q$, 
so these intersections are trivial.
Consider the operators $S=P(P-Q)\inv$ and $S'=-Q(P-Q)\inv$.
We have $\im S = \im P$, $\im S' = \im Q$, and $S+S'=1$, so $\im P+\im Q = H$.
Similarly, consider the operators $T=(P-1)(P-Q)\inv$ and $T'=(1-Q)(P-Q)\inv$.
We have $\im T = \Ker P$, $\im T' = \Ker Q$, and $T+T'=1$, so $\Ker P+\Ker Q = H$.
All four subspaces $\im P$, $\im Q$, $\Ker P$, $\Ker Q$ are closed.
Therefore, both $(\im P, \im Q)$ and $(\Ker P, \Ker Q)$ lie in $\Grr(H)$.

The first equality of \eqref{eq:PQST} implies
$(2S-1)(P-Q) = 2P-(P-Q) = P+Q$.
Note that invertibility of $P+Q$ implies $\im P + \im Q = H$, 
but \textit{does not} imply $\im P \cap \im Q = 0$.
\endproof

\sub{Subspaces of a Hilbert space.}
If $H$ is a Hilbert space, then each closed subspace of $H$ is complemented,
so $\Gr(H)$ is the set of all closed subspaces of $H$.
The map $\im\colon\Proj(H)\to\Gr(H)$ has a natural section
taking a closed subspace $L\subset H$ to the orthogonal projection $P_L$ of $H$ onto $L$.
Applying Proposition \ref{prop:PQST}, we obtain the following result.

\begin{prop}\label{cor:PQST}
  Let $H$ be a Hilbert space. Then the following statements hold:
	\begin{enumerate}\upskip
		\item The pair $(L, M)$ of closed subspaces of $H$ is complementary if and only if $P_L - P_M$ is invertible.
		      If this is the case, then 
	        \begin{equation}\label{eq:PLM}
		         P_{L,M}=P_L(P_L - P_M)\inv.
	        \end{equation}
		\item Let $P\in\Proj(H)$. Then the operator $P+P^*-1$ is invertible, and
		      the orthogonal projection on the image of $P$ is given by the formula
	        \begin{equation}\label{eq:Port}
		         P\ort=P(P+P^*-1)\inv.
	        \end{equation}
	\end{enumerate}	
\end{prop}

\proof
1. \; If $(L,M)\in\Grr(H)$, then also $(L^{\bot},M^{\bot})\in\Grr(H)$.
Applying Proposition \ref{prop:PQST} to the pair of orthogonal projections $P_L$ and $P_M$,
we obtain the first claim of the Corollary.

2. \; $1-P^*$ is the projection on $(\im P)^{\bot}$ along $(\Ker P)^{\bot}$. 
Applying Proposition \ref{prop:PQST} to the pair of projections $P$ and $1-P^*$,
we see that $P+P^*-1 = P-(1-P^*)$ is invertible and
$P(P+P^*-1)\inv$ is the projection on $\im P$ along $(\im P)^{\bot}$.
\endproof

\subsection{The gap topology on $\Gr(H)$} 

For a \textit{Hilbert space} $H$
the map $L\mapsto P_L$ given by the orthogonal projection allows to identify
$\Gr(H)$ with the subspace $\Proj\ort(H)\subset\Proj(H)$ of orthogonal projections in $H$.
The gap topology on $\Gr(H)$ is induced by the norm topology on $\Proj(H)\subset\BB(H)$.

For a \textit{Banach space} $H$ there is no natural section $\Gr(H)\to\Proj(H)$, 
so the definition of the gap topology on $\Gr(H)$ is slightly more complicated in this case. 
Usually the gap topology on $\Gr(H)$ is defined as the topology induced by the gap metric
\begin{multline}
		\hat{\delta}(L_1, L_2) = \max_{i\neq j} 
				\set{ \sup\set{ \dist(u, L_j)\colon u\in L_i, \norm{u}=1 } }, \\
				\hat{\delta}(0, 0) = 0,
				\quad \hat{\delta}(0, L) = 1 \text{ for } L\neq 0. 
\end{multline}
For a Hilbert space $H$ these two definitions of the gap topology coincide.

Proposition \ref{prop:q-top} below gives an equivalent definition of the gap topology 
on the Grassmanian of a Banach space 
in terms of projections, resembling the definition of the gap topology for Hilbert spaces.

The gap topology on $\Gr(H)$ induces the topology on $\Gr^2(H)$ and on its subspace $\Grr(H)$.

\begin{prop}\label{prop:ProjGr}
  Let $H$ be a Banach space. 
	Then the following statements hold:
	\begin{enumerate}\upskip
		\item The map $\im\colon\Proj(H)\to\Gr(H)$ is continuous.
		\item The map $\varphi\colon\Proj(H)\to\Grr(H)$ taking $P\in\Proj(H)$ to $(\im P, \Ker P)\in\Grr(H)$ 
	is a homeomorphism.
		\item $\Grr(H)$ is open in $\Gr^2(H)$.		
	\end{enumerate}
\end{prop}

We first give the proof in the case of a Hilbert space $H$, because it is simpler 
and because we need only this case in the main part of the paper 
as well as in the proofs of all the results below in the context of Hilbert spaces.
After proving the \q{Hilbert case} we give the proof of the general \q{Banach case}.

\medskip

\proof
\noindent 
1. \textit{Suppose first that $H$ is a Hilbert space.}
The map $\im\colon\Proj(H)\to\Gr(H)$ is continuous.
Indeed, it is the composition of the two maps $\Proj(H)\to\Proj\ort(H)$ and $\im\colon\Proj\ort(H)\to\Gr(H)$,
where the first map is given by formula \eqref{eq:Port} and
$\Proj\ort(H)$ is the subspace of $\Proj(H)$ consisting of orthogonal projections. 
The first map is continuous and the second map is an isometry, so their composition is also continuous.

The conjugation by the involution $P\mapsto 1-P$
takes the map $\im\colon \Proj(H)\to\Gr(H)$ 
to the map $\Ker\colon\Proj(H)\to\Gr(H)$, so the second map is also continuous. 
Therefore, $\varphi$ is continuous.
Obviously, $\varphi$ is bijective.

The inverse map $\Grr(H)\to\Proj(H)$ is given by formula \eqref{eq:PLM} and therefore is continuous.
Thus the map $\Proj(H)\to\Grr(H)$ is a homeomorphism.

To prove that $\Grr(H)$ is open in $\Gr^2(H)$, take arbitrary $(L,M)\in\Grr(H)$. 
The operator $P_L-P_M$ is invertible by Corollary \ref{cor:PQST}.
Choose $\varepsilon>0$ such that $2\varepsilon$-neighbourhood of $P_L-P_M$ in $\BB(H)$
consists of invertible operators.
Then for any $L',M'\in\Gr(H)$ such that $\norm{P_L-P_{L'}}<\varepsilon$, $\norm{P_M-P_{M'}}<\varepsilon$
we have 
$$\norm{(P_L-P_M)-(P_{L'}-P_{M'})} \leq \norm{P_L-P_{L'}} + \norm{P_M-P_{M'}} < 2\varepsilon,$$
so $P_{L'}-P_{M'}$ is invertible. 
Applying again Corollary \ref{cor:PQST}, we obtain $(L',M')\in\Grr(H)$.
This completes the proof of the proposition for Hilbert spaces.

\medskip

\noindent 
2. \textit{Let now $H$ be an arbitrary Banach space.}
The continuity of the map $\im\colon\Proj(H)\to\Gr(H)$ follows from the inequality 
$\hat{\delta}(\im P, \im Q) \leq \norm{P-Q}$.
As above, this implies that $\varphi$ is a continuous bijection.
The continuity of the map $\Grr(H)\to\Proj(H)$, $(L,M)\mapsto P_{L,M}$ follows from \cite[Lemma 0.2]{Neu}.
By \cite[Lemma 1 and Theorem 2]{GM59}, $\Grr(H)$ is open in $\Gr^2(H)$.
This completes the proof of the Proposition for Banach spaces.
\endproof
{\sloppy

}

\begin{prop}\label{prop:q-top}
  Let $H$ be a Banach space. 
	Then the gap topology on $\Gr(H)$ coincides with the quotient topology 
	induced by the map $\im\colon\Proj(H)\to\Gr(H)$, $P\mapsto\im P$.
\end{prop}

\proof
The projection $p_1\colon\Gr^2(H)\to \Gr(H)$ onto the first factor is an open continuous map.
By Proposition \ref{prop:ProjGr}, $\Grr(H)$ is open in $\Gr^2(H)$, 
so the restriction of $p_1$ to $\Grr(H)$ is also an open map.
This restriction maps $\Grr(H)$ continuously and surjectively onto $\Gr(H)$.
Therefore, the gap topology on $\Gr(H)$ coincides with the quotient topology 
induced by the map $p_1\colon\Grr(H)\to\Gr(H)$.
To complete the proof, it is sufficient to apply the homeomorphism 
$\varphi\colon\Proj(H)\to\Grr(H)$ from Proposition \ref{prop:ProjGr}.
\endproof

\subsection{Injective maps of Banach spaces}
\label{app:inj}

\begin{prop}\label{prop:Grj}
Let $j\in\BB(H,H')$ be an injective map of Banach spaces.
Denote by $\Gr_j(H)$ the subspace of $\Gr(H)$ consisting of $L$ with $j(L)\in\Gr(H')$.
Then $\Gr_j(H)$ is open in $\Gr(H)$ and the natural inclusion 
$j_*\colon\Gr_j(H) \hookto \Gr(H')$, $L\mapsto j(L)$ is continuous.
\end{prop}

\proof
By Proposition \ref{prop:ProjGr}, $\Gr^M(H)$ is open in $\Gr(H)$.
Thus the statement of the proposition results from the following lemma.

\begin{lem}\label{lem:LMLM}
Let $L\in\Gr_j(H)$, let $M'\in\Gr(H')$ be a complement subspace for $L'=j(L)$, and $M=j\inv(M')$.
Then $L\in\Gr^M(H)\subset\Gr_j(H)$,
and the restriction of $j_*$ to $\Gr^M(H)$ is continuous.
\end{lem}

\prooff{ of the Lemma}
Denote by $Q'$ the projection of $H'$ onto $L'$ along $M'$.
By the Closed Graph Theorem, the bounded linear operator $\restr{j}{L}\colon L\to L'$ is an isomorphism.
Thus the composition $Q = (\restr{j}{L})\inv Q'j$ is a bounded operator on $H$.
Obviously, $Q$ is an idempotent, $\im Q=L$, and $\ker Q=M$.
This implies that $L$ and $M$ are complement subspaces of $H$.

Let $N\in\Gr^M(H)$, $N'=j(N)$.
Then $Q_N = jP_{N,M}(\restr{j}{L})\inv Q'$ is a bounded operator acting on $H'$.
The kernel of $Q_N$ is $M'$ and the restriction of $Q_N^2-Q_N$ to $L'$ vanishes, 
so $Q_N^2=Q_N$ and $Q_N\in\Proj(H')$.
The image of $Q_N$ contains in $N'$ and $N'\cap M' = j(N\cap M)=0$. 
Therefore, $Q_N=P_{N',M'}$, $N' = \im Q_N \in\Gr(H')$, and $N\in\Gr_j(H)$.

By Proposition \ref{prop:ProjGr}, the map $N\mapsto P_{N,M}$ is continuous.
Thus the map $\Gr^M(H)\to\Proj(H')$, $N\mapsto Q_N$ is also continuous.
Composing it with the continuous map $\im\colon\Proj(H')\to\Gr(H')$,
we obtain the continuity of the map $j_*\colon\Gr^M(H)\to\Gr(H')$, $N\mapsto j(N)=\im Q_N$.
This completes the proof of the lemma and of Proposition \ref{prop:Grj}.

\subsection{Closed operators}\label{app:cl}

Let $H$ and $H'$ be Hilbert spaces.
The space $\Cl(H,H')$ of closed linear operators from $H$ to $H'$ is the subspace of $\Gr(H\oplus H')$
consisting of closed subspaces of $H\oplus H'$ which injectively projects on $H$.
An element of $\Cl(H,H')$ can be identified with a linear (not necessarily bounded) operator $A$ acting to $H'$
from (not necessarily closed or dense) subspace $\dom(A)$ of $H$
such that the graph of $A$ is a closed subspace of $H\oplus H'$.

\medskip

All results of this subsection are valid for Banach spaces as well.
\textit{However, in this case the space $\Cl(H,H')$ as we define it (namely, as a the subspace of $\Gr(H\oplus H')$)
does not contain all closed linear operators from $H$ to $H'$, 
but only those whose graphs are complemented subspaces of $H\oplus H'$.}
Nevertheless, families of such operators often arise in applications, so these results can be used for them as well.
For example, for Banach spaces $H$, $H'$ and a linear operator $A$ acting from $\D\subset H$ to $H'$, 
if $\Ker A\subset H$ and $\im A\subset H'$ are closed complemented subspaces,
then the graph of $A$ is a closed complemented subspace of $H\oplus H'$.
In particular, every (not necessarily bounded) Fredholm operator has this property.

\begin{prop}\label{prop:bounded}
Let $H$, $H'$ be Banach spaces.
Then the map $\BB(H,H')\times\Gr(H) \to \Cl(H,H')$ taking $(A,\D)$ to $\restr{A}{\D}$ is continuous.
\end{prop}

\proof
For each $A\in\BB(H,H')$ we define the automorphism $J_A$ of $H\oplus H'$ 
by the formula $J_A(u\oplus u') = u\oplus (u'-Au)$.
Both $A\mapsto J_A$ and $A\mapsto J_A\inv$ are continuous maps from $\BB(H,H')$ to $\BB(H\oplus H')$.
The formula $f(A,Q)= J_A\inv QP_{H,H'}J_A$ defines the continuous map 
$f\colon\BB(H,H')\times\Proj(H) \to \Proj(H\oplus H')$ 
(here $P_{H,H'}$ denotes the projection of $H\oplus H'$ on $H$ along $H'$).
Since $J_A$ takes the graph of $\restr{A}{\D}$ to $\D\oplus 0$ for each $\D\in\Gr(H)$,
$f(A,Q)$ is the projection of $H\oplus H'$ onto the graph of $\restr{A}{\im Q}$.
In other words, we have the commutative diagram
\[
\begin{tikzcd}
  \BB(H,H')\times\Proj(H) \arrow{r}{f} \arrow{d}{\Id\times\im} & \Proj(H\oplus H') \arrow{d}{\im} 
	\\
	\BB(H,H')\times\Gr(H) \arrow{r}{g} & \Gr(H\oplus H')
\end{tikzcd}
\]
where $g$ is the map taking the pair $(A,\D)$ to the graph of $\restr{A}{\D}$.
The top and the right arrows of the diagram are continuous maps, 
while the left arrow is a quotient map by Proposition \ref{prop:q-top}.
Therefore, $g$ is also continuous.
This completes the proof of the proposition.
\endproof

\begin{prop}\label{prop:closed}
Let $W$, $H$, $H'$ be Banach spaces, and let $j\in\BB(W,H)$ be injective. 
Denote by $\Cl_j(W,H')$ the subspace of $\Cl(W,H')$ consisting of operators $A\colon \dom(A)\to H'$
such that the operator $j_*A\colon j(\dom(A))\to H'$, $j_*A = A\cdot j\inv$ lies in $\Cl(H,H')$.
Then the natural inclusion $j_*\colon\Cl_j(W,H') \hookto \Cl(H,H')$ is continuous.
\end{prop}

\proof
Consider the following commutative diagram:
\[
\begin{tikzcd}
  \Cl(W,H') \arrow[hookleftarrow]{r} \arrow[hookrightarrow]{d} &
	  \Cl_j(W,H') \arrow{r}{j_*} \arrow[hookrightarrow]{d} &
	  \Cl(H,H') \arrow[hookrightarrow]{d} \\
	\Gr(W \oplus H') \arrow[hookleftarrow]{r}  &
	  \Gr_j(W \oplus H') \arrow{r}{j_*} & \Gr(H \oplus H')
\end{tikzcd}
\]
The spaces above are just subspaces on the spaces below,
and $\Cl_j(W,H') = \Cl(W,H') \cap \Gr_j(W \oplus H')$.
By Proposition \ref{prop:Grj}, the map $j_*\colon \Gr_j(W \oplus H') \to \Gr(H \oplus H')$ is continuous.
So the restriction of $j_*$ to $\Cl_j(W,H') \subset \Gr_j(W \oplus H')$ is also continuous.
This completes the proof of the proposition.
\endproof

\subsection{Differential and pseudo-differential operators}
\label{sec:app-diff}

The results of the previous subsection can be used for
differential and pseudo-differen\-tial operators acting between sections of vector bundles over $M$.
To achieve continuity of the corresponding families of closed operators,
the relevant topology on the space of differential operators 
will be the $C^0_b$-topology on their coefficients.

\sub{General framework.}
Let $M$ be a smooth Riemannian manifold and $E$, $E'$ be smooth Hermitian vector bundles over $M$.
For an integer $\d\geq 1$,
we denote by $\Opd(E,E')$ the set of all pairs $(A,\D)$ such that
\begin{itemize}\upskip
	\item\upskip $A$ is a bounded operator from $H^d(E)$ to $L^2(E')$, 
	\item\upskip $\D$ is a closed subspace of $H^d(E)$, and
	\item\upskip the restriction $\restr{A}{\D}$ of $A$ to the domain $\D$ is closed as an operator from $L^2(E)$ to $L^2(E')$.
\end{itemize}\uppskip
We equip $\Opd(E,E')$ with the topology induced by the inclusion
$$\Opd(E,E') \hookto  \BB(H^d(E), L^2(E')) \times \Gr(H^d(E)).$$
Here $L^2(E)$ is the Hilbert space of $L^2$-sections of $E$
and $H^d(E)$ is the $d$-th order Sobolev space of sections of $E$.

\begin{prop}\label{prop:diff}
The  map $\Opd(E,E') \to \Cl(L^2(E),L^2(E'))$ 
taking $(A,\D)$ to $\restr{A}{\D}$
is continuous.
\end{prop}

\proof
Take $W=H^d(E)$, $H=L^2(E)$, $H'=L^2(E')$,
and let $j$ be the natural embedding $W \hookto H$.
By Proposition \ref{prop:bounded}, the map 
$\Opd(E,E') \subset \BB(W,H')\times\Gr(W) \to \Cl(W,H')$ is continuous.
By definition of $\Opd(E,E')$, the image of this map is contained in $\Cl_j(W,H')$.
By Proposition \ref{prop:closed},
the map $j_*\colon\Cl_j(W,H') \to \Cl(H,H')$ is continuous.
Combining all this, we obtain the continuity of the map $\Opd(E,E') \to \Cl(H,H')$.
\endproof

This general result can be applied to differential or pseudo-differential operators $A$ of order $d$ 
with domains $\D$ given by boundary conditions.
We show below how Proposition \ref{prop:diff} can be applied to 
boundary value problems for first order differential operators,
in particular local boundary value problems.
We omit discussion of higher order operators 
because boundary conditions are slightly more complicated in that case;
however, Proposition \ref{prop:diff} works for higher order operators as well.

\sub{Boundary value problems for first order operators.}
Suppose now that $M$ is a compact manifold.
Denote by $\EY$ the restriction of $E$ to the boundary $\pM$.

Let $A\in\BB(H^1(E), L^2(E'))$.
In particular, $A$ can be a first order differential operator with continuous coefficients.
For a closed subspace $\L$ of $H^{1/2}(\EY)$ 
we denote by $A_{\L}$ the operator $A$ with the domain
$$\dom\br{A_{\L}} = \set{ u\in H^1(E)\colon \tau(u) \in\L},$$
where $\tau\colon H^1(E) \to H^{1/2}(\EY)$ is the trace map 
extending by continuity the restriction map $C^{\infty}(E) \to C^{\infty}(\EY)$, $u\mapsto\restr{u}{\pM}$.

Let $\Opt(E,E')$ denotes the subspace of $\BB(H^1(E), L^2(E')) \times \Gr(H^{1/2}(\EY))$
consisting of pairs $(A,\L)$ such that the operator $A_{\L}$ is closed.

\begin{prop}\label{prop:trace}
The map 
$$\Opt(E,E')\to\Cl(L^2(E),L^2(E')), \quad (A,\L) \mapsto A_{\L}$$
is continuous.
\end{prop}

\proof
The inverse image $\tau\inv(\L)$ is a closed subspace of $H^1(E)$.
Since $\tau$ is bounded and surjective, the map
$$\tau^*\colon \Gr(H^{1/2}(\EY)) \to \Gr(H^1(E)), \quad \L \mapsto \tau\inv(\L),$$
is continuous.
Hence the map $\Opt(E,E')\to\Op1(E,E')$ taking $(A,\L)$ to $(A,\tau\inv(\L))$ is also continuous.
It remains to apply Proposition \ref{prop:diff}.
\endproof

\sub{Local boundary value problems for first order operators.}
Denote by $\Ellz(E,E')$ the set of all pairs $(A,L)$ such that
\begin{itemize}\upskip
	\item $A$ is a first order elliptic differential operator with smooth coefficients 
acting from sections of $E$ to sections of $E'$, and
	\item $L$ is a smooth subbundle of $\EY$ satisfying Shapiro-Lopatinskii condition \eqref{eq:Lell-gen}.
\end{itemize}\upskip
Equip $\Ellz(E,E')$ with the $C^0$-topology on coefficients of operators 
and the $C^1$-topology on boundary conditions,
that is, the topology induced by the inclusion
$$\Ellz(E,E') \hookto \BB(H^1(E),L^2(E')) \times C^1(\Gr(\EY)).$$
Here $\Gr(\EY)$ denotes the smooth vector bundle over $\pM$ 
whose fiber over $x\in \pM$ is the Grassmanian $\Gr(E_x)$,
and sections of $\Gr(\EY)$ are identified with subbundles of $\EY$.

The Sobolev space $H^{1/2}(L)$ can be naturally identified with the closed subspace of $H^{1/2}(\EY)$
via the map $H^{1/2}(L)\ni u \mapsto u\oplus 0 \in H^{1/2}(L)\oplus H^{1/2}(L^{\bot}) = H^{1/2}(\EY)$.
This allows to associate with a pair $(A,L)\in\Ellz(E,E')$
the unbounded operator $A_L$ acting as $A$ on the domain
$$ \dom\br{A_L} = \set{ u\in H^1(E)\colon \tau(u) \in H^{1/2}(L)}.$$
By the classical theory of elliptic operators, $A_L$ is closed for every $(A,L)\in\Ellz(E,E')$.
See, for example, Proposition \ref{prop:AL1}, where it is proven for self-adjoint operators.
Closedness of a non-self-adjoint $A_L$ can be proven along the same lines, 
or can be obtained directly from Proposition \ref{prop:AL1}
by replacing a pair $(A,L)$ with the pair $(A',L')\in\Ellz(E\oplus E')$, where
$A'=\smatr{ 0 & A^t \\ A & 0}$ and $L' = L\oplus (\sigma_A(n)L)^{\bot} \subset \EY\oplus\EY'$.

\begin{prop}\label{prop:AL2-app}
The natural inclusion $\Ellz(E,E') \hookto \Cl\br{L^2(E),L^2(E')}$, $(A,L) \mapsto A_L$ 
is continuous.
\end{prop}

\proof
It is an immediate corollary of the following lemma applied to $N=\pM$ and $F=\EY$
and of Proposition \ref{prop:trace}.
\endproof

\begin{lem}\label{lem:local}
Let $F$ be a smooth Hermitian vector bundle over a smooth closed Riemannian manifold $N$.
Then the map
\begin{equation}\label{eq:LH12}
  C^{\infty,1}(\Gr(F)) \to \Gr(H^{1/2}(F)),
\end{equation}
taking a smooth subbundle $L$ of $F$ to $H^{1/2}(L) \subset H^{1/2}(F)$,
is continuous.
Here $C^{\infty,1}(\Gr(F))$ denotes the space of smooth sections of $\Gr(F)$ with the $C^{1}$-topology,
that is, the topology induced by the embedding $C^{\infty}(\Gr(F)) \hookto C^{1}(\Gr(F))$.
\end{lem}

\proof
Operator of multiplication by a $C^1$-function $N\to\CC$ is a bounded operator on $H^s(N)$ for every $s\in[0,1]$.
In particular, it is bounded as an operator acting on $H^{1/2}(N)$,
and the correspondent inclusion $C^1(N) \hookto \BB\br{H^{1/2}(N)}$ is continuous.
Passing to bundles, we obtain the natural continuous inclusion
\begin{equation}\label{eq:C1H}
  C^1(\End(F)) \hookto \BB\br{H^{1/2}(F)}.
\end{equation}

The smooth map $P\colon\Gr(\CC^{n})\to\End(\CC^{n})$, $V\mapsto P_V$,
induces the continuous map
\[
P_* \colon C^{1}(\Gr(F)) \hookto C^{1}(\End(F)),
\]
which carries a subbundle $L$ of $F$ to the orthogonal projection $P_*L$ of $F$ onto $L$.
Composing it with the continuous inclusion \eqref{eq:C1H},
we obtain the continuous map
\[
Q\colon C^{1}(\Gr(F)) \hookto \BB\br{H^{1/2}(F)}.
\]
For each smooth subbundle $L$ of $F$ the bounded operator $Q(L)$ 
is an idempotent with the image $H^{1/2}(L)$.
By Proposition \ref{prop:ProjGr}(1), the map
$$\im\colon \Proj\br{H^{1/2}(F)} \to \Gr\br{H^{1/2}(F)},$$
is continuous.
Composing it with $Q$, we obtain continuity of \eqref{eq:LH12}.
This completes the proof of the lemma.
\endproof

\bigskip
\bigskip
\noindent
\textit{Department of Mathematics, Technion -- Israel Institute of Technology, Haifa, Israel} \\
\textit{marina.p@campus.technion.ac.il} 


\begin{thebibliography}{99}

\bibitem{BLP-04}
 B. Booss-Bavnbek, M. Lesch, and J. Phillips.
 Unbounded Fredholm operators and spectral flow.
 Canadian Journal of Mathematics \textbf{57}(2), 225--250 (2005);
 arXiv:math/0108014 [math.FA].

\bibitem{BBW-93}
 B. Booss-Bavnbek, K. P. Wojciechowski.
 Elliptic boundary problems for Dirac operators.
 Birkhauser, 1993.

\bibitem{BBZ-14}
B. Booss-Bavnbek, Chaofeng Zhu.
The Maslov index in symplectic Banach spaces.
Memoirs of the American Mathematical Society, \textbf{252}(1201), 1--118 (2018);
arXiv:1406.0569 [math.SG].

\bibitem{CL}
H. O. Cordes, J. P. Labrousse.
The invariance of the index in the metric space of closed operators.
Journal of Mathematics and Mechanics \textbf{12}(5), 693--719 (1963).

\bibitem{GM59}
I. Z. Gohberg, A. S. Markus. 
Two theorems on the gap between subspaces of a Banach space. 
Uspehi Mat. Nauk \textbf{14}:5(89), 135--140 (1959) (in Russian).

\bibitem{GL}
A. Gorokhovsky, M. Lesch.
On the spectral flow for Dirac operators with local boundary conditions.
International Mathematics Research Notices, \textbf{17}, 8036-8051 (2015);
arXiv:1310.0210 [math.AP].

\bibitem{Hor}
L. H\"ormander. 
The analysis of linear partial differential operators III. 
Grundlehren der Mathematischen Wissenschaften, vol. 274, Springer-Verlag, Berlin--Heidelberg--New York, 1985.

\bibitem{Jo}
M. Joachim.
Unbounded Fredholm operators and K-theory.
In: High-dimen\-sio\-nal Manifold Topology, World Sci. Publ., River Edge, NJ, 177--199 (2003).

\bibitem{Kato}
T. Kato. 
Perturbation Theory for Linear Operators. 
A Series of Comprehensive Studies in Mathematics 132 (1980).

\bibitem{KN}
M. I. Katsnelson, V. E. Nazaikinskii.
The Aharonov-Bohm effect for massless Dirac fermions and the spectral flow of Dirac type operators with classical boundary conditions.
Theoretical and Mathematical Physics \textbf{172}(3), 1263--1277 (2012);
arXiv:1204.2276 [math.AP].

\bibitem{KL}
P. Kirk, M. Lesch.
The eta-invariant, Maslov index, and spectral flow for Dirac-type operators on manifolds with boundary.
Forum Math \textbf{16}, 553--629 (2004);
arXiv:math/0012123 [math.DG].

\bibitem{Lax}
P. D. Lax, R. S. Phillips. 
Local boundary conditions for dissipative symmetric linear differential operators. 
Communications on Pure and Applied Mathematics, \textbf{13}(3), 427--455 (1960).

\bibitem{L-04}
M. Lesch.
The uniqueness of the spectral flow on spaces of unbounded self-adjoint Fredholm operators.
In: Spectral geometry of manifolds with boundary and decomposition of manifolds
(B.~Booss-Bavnbek, G.~Grubb, and K.P. Wojciechowski, eds.),
AMS Contemporary Math Proceedings \textbf{366}, 193--224  (2005);
arXiv:math/0401411 [math.FA].

\bibitem{Neu}
G. Neubauer.
Homotopy properties of semi-Fredholm operators in Banach spaces. 
Mathematische Annalen \textbf{176}(4), 273--301 (1968). 

\bibitem{Ph}
J. Phillips.
Self-adjoint Fredholm operators and spectral flow.
Canadian Mathematical Bulletin \textbf{39}(4), 460--467 (1996).

\bibitem{Pr}
M. Prokhorova.
The spectral flow for Dirac operators on compact planar domains with local boundary conditions.
Communications in Mathematical Physics \textbf{322}(2), 385--414 (2013);
arXiv:1108.0806 [math-ph].

\bibitem{Pr14}
M. Prokhorova.
The spectral flow for elliptic operators on compact surfaces.
Current problems in mathematics, computer science and natural science knowledge.
Proceedings of International scientific conference dedicated to the 155th anniversary
of I.V. Meshcherskiy (September 2014, Koryazhma, Russia),
78--80 (2014).

\bibitem{Pr17}
M. Prokhorova.
The spectral flow for local boundary value problems on compact surfaces. 
arXiv:1703.06105v1 [math.AP] (2017).

\bibitem{Pr18}
M. Prokhorova.
Self-adjoint local boundary problems on compact surfaces. II. Family index.
arXiv:1809.04353 [math-ph] (2018).

\end{thebibliography}
\end{document}